\input{epsf}
\newcommand{\et}{\hspace{-0.08in}{\bf .}\hspace{0.1in}}
\newcommand{\BOX}{\hbox {$\sqcap$ \kern -1em $\sqcup$}}
\newcommand{\qed}{\hskip 2em \hbox{\BOX} \vskip 2ex}

\newcommand{\To}{\Rightarrow}
\newcommand{\Set}{{\rm Set}}
\newcommand{\Cat}{{\rm Cat}}
\newcommand{\Top}{{\rm Top}}
\newcommand{\TopCat}{{\rm TopCat}}
\newcommand{\Grp}{{\rm Grp}}
\newcommand{\Diff}{{\rm Diff}}
\newcommand{\DiffCat}{{\rm DiffCat}}

\renewcommand{\to}{\rightarrow}
\newcommand{\tensor}{\otimes}
\newcommand{\maps}{\colon}
\newcommand{\inv}{{\rm inv}}

\newcommand{\iso}{\cong}

\newcommand{\id}{{\rm id}}

\newcommand{\wg}{{\rm W2G }} 
\newcommand{\cg}{{\rm C2G }}

\newcommand{\xb}{\bar{x}}

\newcommand{\ten}{\otimes }
\newcommand{\Mon}{{\rm Mon   Cat}}
\newcommand{\imp}{{\rm Imp}}

\newtheorem{thm}{Theorem}    

\newtheorem{prop}[thm]{Proposition}
\newtheorem{defn}[thm]{Definition}

        \newcommand{\be}{\begin{equation}}
        \newcommand{\ee}{\end{equation}}
        \newcommand{\ba}{\begin{eqnarray}}
        \newcommand{\ea}{\end{eqnarray}}
        \newcommand{\ban}{\begin{eqnarray*}}
        \newcommand{\ean}{\end{eqnarray*}}
        \newcommand{\barr}{\begin{array}}
        \newcommand{\earr}{\end{array}}

\documentclass{article}
\usepackage{amsfonts,amssymb}
\usepackage{diagram}
\usepackage[all]{xy}
   \title{Remarks on 2-Groups}
      \author{Aaron D.\ Lauda \\
      Department of Physics,  University of California\\
      Riverside, California 92521 \\
      USA \\
         \\
      email: lauda@math.ucr.edu \\}

\begin{document}
\bibliographystyle{plain}
\maketitle

\begin{abstract}

\noindent A 2-group is a `categorified' version of a group, in
which the underlying set $G$ has been replaced by a category and
the multiplication map $m \maps G \times G \to G$ has been
replaced by a functor.  A number of precise definitions of this
notion have already been explored, but a full treatment of their
relationships is difficult to extract from the literature.  Here
we describe the relation between two of the most important
versions of this notion, which we call `weak' and `coherent'
2-groups.  A weak 2-group is a weak monoidal category in which
every morphism has an inverse and every object $x$ has a `weak
inverse': an object $y$ such that $x \tensor y \iso 1 \iso
y \tensor x$. A coherent 2-group is a weak 2-group in which every
object $x$ is equipped with a {\it specified} weak inverse $\xb$
and isomorphisms $i_x \maps 1 \to x \tensor \xb$, $e_x \maps \xb
\tensor x \to 1$ forming an adjunction.  We define 2-categories of
weak and coherent 2-groups and construct an `improvement'
2-functor which turns weak 2-groups into coherent ones; using
this one can show that these 2-categories are biequivalent.  We also
internalize the concept of a coherent 2-group.  This gives a way
of defining topological 2-groups, Lie 2-groups, and the like.
\end{abstract}

\section{Introduction}

Group theory has proven to be a powerful tool not only in mathematics, but
also in physics, chemistry and other sciences.  In recent times it has
become evident that in many contexts where we are tempted to use groups,
it is actually more natural to use a richer sort of structure, namely a
kind of `higher-dimensional' group.  One might also call this a
`categorified' group, since the underlying set $G$ of a traditional
group has been replaced by a {\it category} and the multiplication function
$m \maps G \times G \to G$ has been replaced by a {\it functor}.  To
hint at a sequence of further generalizations where we use $n$-categories
and $n$-functors, we call this sort of thing a `2-group'.

There are various different ways to make the concept of 2-group more
precise.  Some can already be found in the mathematical literature, but
unfortunately they often remain implicit in work that focuses on more
general concepts.  Indeed, while the basic facts about 2-groups are
familiar to most experts in category theory, it is impossible for
beginners to find a unified presentation of this material with all
the details provided.  The present paper tries to start filling this gap.

Whenever one categorifies a mathematical concept, there are some choices
involved.  For example, one might define a 2-group simply to be a
category $G$ equipped with functors describing multiplication, inverses
and the identity, satisfying the usual group axioms `on the nose' ---
that is, as equations between functors.  We call this a `strict'
2-group.  Strict 2-groups have been applied in a variety of contexts,
including homotopy theory \cite{Brown,BS}, topological quantum field
theory \cite{Yetter}, and gauge theory \cite{Baez}.  Part of the charm
of strict 2-groups is that they can be defined in a large number of
equivalent ways, including:
\begin{itemize}
    \item a strict monoidal category in which all objects
     and morphisms are invertible,
    \item a strict 2-category with one object in which all 1-morphisms and
    2-morphisms are invertible,
    \item a group object in $\Cat$ (also called a `categorical group')
    \item a category object in $\Grp$,
    \item a crossed module.
\end{itemize}
There is an excellent review article by Forrester-Barker that
explains most of these notions and why they are equivalent \cite{FB}.

However, as the notion of a group takes on this higher dimensional form,
the exact definition most suited for a given task becomes less obvious.
For instance, rather than imposing the group axioms as equational laws,
we could instead require that they hold up to specified isomorphisms
satisfying laws of their own.  This leads to the concept of a
`coherent 2-group'.

For example, given objects $x,y,z$ in a strict 2-group we have
\[ (x \tensor y) \tensor z = x \tensor (y \tensor z)  \]
where we write multiplication as $\tensor$.   In a coherent 2-group,
we instead specify an isomorphism called the `associator':
\[
a_{x,y,z} \maps
\xymatrix@1{(x \ten y) \ten z \; \ar[r]^<<<<<{\sim} & \; x \ten (y \ten z)} .
\]
Similarly, we replace the left and right unit laws
\[        1 \tensor x = x, \qquad x \tensor 1 = x \]
by isomorphisms
\[
   \ell_x \maps \xymatrix@1{1 \tensor x \; \ar[r]^<<<<<{\sim} & \; x}
, \qquad
     r_x \maps \xymatrix@1{x \tensor 1\; \ar[r]^<<<<<{\sim} &\; x}
\]
and replace the equations
\[ x \tensor x^{-1} = 1, \qquad x^{-1} \tensor x = 1 \]
by isomorphisms called the `unit' and `counit'.

Next, in order to manipulate these isomorphisms with some of the same
facility as with equations, we require that they satisfy conditions
known as `coherence laws'.  The coherence laws for the associator and
the left and right unit laws were developed by Mac Lane \cite{MacLane}
in his groundbreaking work on monoidal categories, while those for the
unit and counit are familiar from the definition of an adjunction in a
monoidal category \cite{JS}.  Putting these ideas together, one obtains
Ulbrich and Laplaza's definition of a `category with group structure'
\cite{Laplaza,Ulbrich}.  Finally, a `coherent 2-group' is a category $G$
with group structure in which all morphisms are invertible.  This last
condition ensures that there is a covariant functor
\[      \inv \maps G \to G   \]
sending each object $x \in G$ to its weak inverse $\xb$; otherwise
there will only be a contravariant functor of this sort.

In this paper we compare this sort of 2-group to a simpler sort, which
we call a `weak 2-group'.  This is a weak monoidal category in which
every morphism has an inverse and every object $x$ has a `weak inverse':
an object $y$ such that $y \tensor x \iso 1$ and $x \tensor
y \iso 1$.  Note that in this definition, we do not {\it specify} the
weak inverse $y$ or the isomorphisms from $y \tensor x$ and $x \tensor
y$ to $1$, nor do we impose any coherence laws upon them.  Instead, we
merely demand that they {\it exist}.  Nonetheless, it turns out that any
weak 2-group can be improved to become a coherent one!  While this
follows immediately from a theorem of Laplaza \cite{Laplaza}, it
seems worthwhile to give an expository account here, and to formalize
this process as a 2-functor
\[ \imp \maps \wg \to \cg \]
where $\wg$ and $\cg$ are suitable 2-categories of weak and coherent
2-groups, respectively.

To do this, we start in Section \ref{weaksection} by defining weak
2-groups and the 2-category $\wg$. In Section \ref{coherentsection} we
define coherent 2-groups and the 2-category $\cg$.  In Section
\ref{internalizationsection} we show that the concept of `coherent
2-group object' can be defined in any 2-category with finite
products. This allows us to define notions such as `coherent topological
2-group', `coherent Lie 2-group' and the like.  While this may seem a
bit of a digression, it serves as an excellent excuse to introduce the
technique of string diagrams \cite{Street}, which turn out to be crucial
for constructing the 2-functor $\imp \maps \wg \to \cg$.  We
construct this 2-functor in Section \ref{improvementsection}.
Together with the forgetful 2-functor ${\rm F} \maps \cg \to \wg$,
this sets up a `biequivalence' between $\wg$ and $\cg$.

In other words, the 2-category of weak 2-groups and the 2-category
of coherent 2-groups are `the same' in a suitably weakened sense.
Thus there is not really too much difference between weak and coherent
2-groups: we can freely translate theorems about
one into theorems about the other using the 2-functors
$\imp \maps \wg \to \cg$ and ${\rm F} \maps \cg \to \wg$.

{\bf Note}: in all that follows, we write the composite of morphisms
$f \maps x \to y$ and $g \maps y \to z$ as $fg \maps x \to z$.

\section{Weak 2-groups} \label{weaksection}

Before we define a weak 2-group, recall that a
{\bf weak monoidal category} consists of:
\begin{description}
    \item[(i)] a category $M$,
    \item[(ii)] a functor $m \maps M \times M \to M$, where we
    write $m(x,y)=x \ten y$ and $m(f,g)=f \ten g$ for objects $x, y,
    \in M$ and morphisms $f, g$ in $M$,
    \item[(iii)] an `identity object' $1 \in M$,
    \item[(iv)] natural isomorphisms
     \[ a_{x,y,z} \maps (x \ten y) \ten z \to x \ten (y \ten z), \]
     \[ \ell_x \maps 1 \ten x \to x , \]
     \[ r_x \maps x \ten 1 \to x,  \]
\end{description}
such that the following diagrams commute for all objects $x,y,z,w \in M$:
\[
 \xy
    (0,20)*{(x \ten y)(z \ten w)}="1";
    (40,0)*{x \ten (y \ten (z \ten w))}="2";
    (25,-20)*{ \quad x \ten ((y \ten z) \ten w)}="3";
    (-25,-20)*{(x \ten (y \ten z)) \ten w}="4";
    (-40,0)*{((x \ten y) \ten z) \ten w}="5";
        {\ar^{a_{x,y,z \ten w}}     "1";"2"}
        {\ar_{1_x \ten a _{y,z,w}}  "3";"2"}
        {\ar^{a _{x,y \ten z,w}}    "4";"3"}
        {\ar_{a _{x,y,z} \ten 1_w}  "5";"4"}
        {\ar^{a _{x \ten y,z,w}}    "5";"1"}
 \endxy
 \\
\]
\vskip 1em
\[
 \xymatrix{
    (x \ten 1) \ten y
        \ar[rr]^{a_{x,1,y}}
        \ar[dr]_{r_x \ten y}
     &&  x \ten (1 \ten y)
        \ar[dl]^{x \ten \ell_y } \\
    & x \ten y   } \\
\]
A {\bf strict monoidal category} is the special case where
$a_{x,y,z},\ell_x,r_x$ are all identity morphisms.  In this
case we have
\[ (x \tensor y) \tensor z = x \tensor (y \tensor z) , \]
\[        1 \tensor x = x, \qquad x \tensor 1 = x .\]
As mentioned in the Introduction, a {\bf strict 2-group} is a strict
monoidal category where every morphism is invertible and every
object $x$ has an inverse $x^{-1}$, meaning that
\[      x \tensor x^{-1} = 1, \qquad x^{-1} \tensor x = 1  .\]

Following the principle that it is wrong to impose equations between
objects in a category, we can instead start with a weak monoidal
category and require that every object has a `weak' inverse.  With these
changes we obtain the definition of `weak 2-group':

\begin{defn} \label{weakinv} \et
If $x$ is an object in a weak monoidal category, a
{\bf weak inverse} for $x$ is an object $y$ such that
$x \tensor y \iso 1$ and $y \tensor x \iso 1$.  If $x$
has a weak inverse, we call it {\bf weakly invertible}.
\end{defn}

\begin{defn} \label{weak} \et
A {\bf weak 2-group} is a weak monoidal category where all
objects are weakly invertible and all morphisms are invertible.
\end{defn}

Weak 2-groups are the objects of a 2-category $\wg$; now let us
describe the morphisms and 2-morphisms in this 2-category.  Notice that
the only {\it structure} in a weak 2-group is that of its underlying
weak monoidal category; the invertibility conditions on objects and
morphisms are only {\it properties}.  With this in mind, it is natural
to define a morphism between weak 2-groups to be a weak monoidal
functor.  Recall that a {\bf weak monoidal functor} $F \maps C \to C'$
between monoidal categories $C$ and $C'$ consists of:
\begin{description}
\item[(i)] a functor $F \maps C \to C'$,
\item[ii)] a natural isomorphism $F_{2} \maps F(x) \otimes F(y)
\to F(x \otimes y)$, where for brevity we suppress the subscripts
indicating the dependence of this isomorphism on $x$ and $y$,
\item[(iii)] an isomorphism $F_{0} \maps 1' \to F(1)$,
where $1$ is the unit object of $C$ and $1'$ is the unit object of $C'$,
\end{description}
such that the following diagrams commute for all objects $x,y,z \in C$:
\[
 \xymatrix@!C{
    (F(x) \ten F(y)) \ten F(z)
      \ar[r]^>>>>>>>{F_{2} \ten 1}
      \ar[d]_{a_{F(x), F(y), F(z)}}
 & F(x \ten y) \ten F(z)
      \ar[r]^{F_{2}}
 & F((x \ten y) \ten z)
      \ar[d]^{F(a_{x,y,z})}   \\
    F(x) \ten (F(y) \ten F(z))
      \ar[r]^>>>>>>>{1 \ten F_{2}}
 & F(x) \ten F(y \ten z)
      \ar[r]^{F_{2}}
 & F(x \ten (y \ten z))  .
    }
\]
\vskip 1em
\[
 \xymatrix{
    1' \ten F(x)
      \ar[r]^{\ell_{F(x)}}
      \ar[d]_{F_{0} \ten 1}
 &  F(x)  \\
    F(1) \ten F(x)
      \ar[r]^{F_{2}}
 &  F(1 \ten x)
      \ar[u]_{F(\ell_{x})}
 }
\]
\vskip 1em
\[
 \xymatrix{
    F(x) \ten 1'
      \ar[r]^{r_{F(x)}}
      \ar[d]_{1 \ten F_{0}}
 &  F(x)    \\
    F(x) \ten F(1)
      \ar[r]^{F_{2}}
 &  F(x \ten 1)  .
      \ar[u]_{F(r_x)}
 }
\]

A weak monoidal functor preserves all the structure of a weak monoidal
category up to specified isomorphisms.  Moreover, if $C$
and $C'$ are weak 2-groups, a weak monoidal functor $F \maps C \to C'$
also preserves weak inverses:

\begin{prop} \et
If $F \maps C \to C'$ is a monoidal functor between weak 2-groups
$C$ and $C'$ and $y \in C$ is a weak inverse of $x \in C$, then
$F(y)$ is a weak inverse of $F(x)$ in $C'$.
\end{prop}

\noindent {\bf Proof. }
Since $y$ is a weak inverse of $x$, there
must exist isomorphisms $\gamma \maps x \ten y \to 1$ and
$\xi \maps y \ten x \to 1$. The proposition is then established by
composing the following isomorphisms:
\[
 \xymatrix{
    F(y) \ten F(x)
      \ar[r]^<<<<<{\sim}
      \ar[d]_{F_{2}}
 &  1       \\
    F(y \ten x)
      \ar[r]_{F(\xi)}
 &  F(1) ,
      \ar[u]_{F_{0}}
  }
 \qquad \qquad
 \xymatrix{
    F(x) \ten F(y)
      \ar[r]^<<<<<{\sim}
      \ar[d]_{F_{2}}
 &  1       \\
    F(x \ten y)
      \ar[r]_{F(\gamma)}
 &  F(1)
      \ar[u]_{F_{0}}
 }
\]
\hskip 30em \qed

We thus make the following definition:

\begin{defn} \et
A {\bf homomorphism} $F \maps C \to C'$ between weak 2-groups
is a weak monoidal functor.
\end{defn}
\noindent
The composite of weak monoidal functors is again a
monoidal functor \cite{moncat}, and composition satisfies
associativity and the unit laws.  Thus, 2-groups and
the homomorphisms between them form a category.

Although no direct counterpart can be found in traditional group theory,
it is natural in this categorified context to also consider
`2-homomorphisms' between homomorphisms.  Since a homomorphism between
weak 2-groups is just a weak monoidal functor, it makes sense to define
2-homomorphisms to be weak monoidal natural transformations.  Recall
that if $F,G \maps C \to C'$ are are weak monoidal functors,
then a {\bf weak monoidal natural
transformation} $\theta \maps F \To G$ is a natural transformation such
that these diagrams commute for all $x,y \in C$:
\[
 \xymatrix{
    F(x) \ten F(y)
     \ar[rr]^{\theta_{x} \ten \theta_{y}}
     \ar[d]_{F_{2}}
 &&  G(x) \ten G(y)
      \ar[d]^{G_{2}}     \\
    F(x \ten y)
      \ar[rr]^{\theta_{x \ten y}}
 &&  G(x \ten y)
 }
\]

\[
\\
\\
 \xymatrix{
    1
      \ar[d]_{F_{0}}
      \ar[dr]^{G_{0}}  \\
    F(1)
      \ar[r]^{\theta_{1}}
 &  G(1)
 }
\]
commute.  Thus we make the following definitions:

\begin{defn} \et
A {\bf 2-homomorphism} $\theta \maps F \To G$ between
homomorphisms \hfill \break $F,G \maps C \to C'$ of weak 2-groups
is a weak monoidal natural transformation.
\end{defn}

\begin{defn} \et Let {\bf W2G} be the 2-category consisting of weak
2-groups, homomorphisms between these, and 2-homomorphisms between
those.  \end{defn}

There is a 2-category $\Mon$ with monoidal categories as objects, weak
monoidal functors as 1-morphisms and weak monoidal natural
transformations as 2-morphisms \cite{moncat}.  Thus $\wg$ a 2-category,
since it is a full and 2-full sub-2-category of $\Mon$.

\section{Coherent 2-groups} \label{coherentsection}

In this section we explore another notion of 2-group.
Rather than requiring that objects be weakly invertible,
we will require that every object be equipped with a
specified adjunction.  Recall that an {\bf adjunction} is a
quadruple $(x, \xb, i_x, e_x)$ where $i_x \maps$ $\xymatrix@1{1
\ar[r]^<<<<<{\sim} & x \ten \xb}$ (called the {\bf unit})
and $e_x \maps$ $\xymatrix@1{\xb \ten x \ar[r]^<<<<<{\sim} & 1}$
(called the {\bf counit}) are morphisms such that the
following diagrams
\[
 \xymatrix@!C{
     1 \ten x \ar[r]^<<<<<<{i_x \ten 1} \ar[d]_{\ell_x}
     & (x \ten \xb) \ten x \ar[r]^{a_{x, \bar{x},x}}
     & x\ten( \bar{x}\ten x) \ar[d]^{1 \ten e_x} \\
     x \ar[rr]_{r^{-1}_x}
      && x \ten 1  }
\]

\[
 \xymatrix@!C{
    \xb \ten 1
        \ar[r]^<<<<<<{1 \ten i_x}
        \ar[d]_{r_{\xb}}
      & \xb \ten (x \ten \xb)
        \ar[r]^{a^{-1}_{\xb, x, \xb}}
      & (\xb \ten x)\ten \xb
        \ar[d]^{e_x \ten 1} \\
    \xb
        \ar[rr]_{\ell^{-1}_{\xb}}
      && 1 \ten \xb   }
\]
commute.  For reasons that will become apparent in the sections to come we
refer to these diagrams as the first and second
{\bf zig-zag identities}, respectively.

An adjunction $(x,\xb,i_x,e_x)$ for which the unit and counit are
invertible is called an {\bf adjoint equivalence}.  In this case $x$ and
$\xb$ are weak inverses.  Thus, specifying an adjoint equivalence for $x$
ensures that $\xb$ is weakly invertible --- but it does so by providing
$x$ with extra {\it structure}, rather than merely asserting a {\it
property} of $x$.  We now make the following definition:

\begin{defn}  \label{coherent} \et
A { \bf coherent 2-group } is a weak monoidal category $C$ in which
every morphism is invertible and every object $x \in C$ is equipped with
an adjoint equivalence $(x,\bar{x},i_{x},e_{x})$.
\end{defn}

\noindent As noted in the Introduction, a coherent 2-group is
the same as a category with group structure \cite{Laplaza,Ulbrich}
in which all morphisms are invertible.   It is also the same as
an `autonomous monoidal category' \cite{JS} with all morphisms
invertible, or a `bigroupoid' \cite{HKK} with one object.

As we did with weak 2-groups, we can define homomorphisms between
coherent 2-groups.  As in the weak 2-group case we can begin with a weak
monoidal functor, but now we must consider what additional structure
this must have to preserve each adjoint equivalence $(x,\xb,i_x,e_x)$,
at least up to a specified isomorphism.  At first it may seem that
an additional structural map is required.
That is, if $F \maps C \to C'$ is a weak monoidal functor
it may seem that we must include a natural isomorphism
\[ F_{-1} \maps \overline{F(x)} \to F(\xb) \]
relating the inverse of the image of $x$ to the image of
the inverse $\xb$.  We shall show this is not the case:
$F_{-1}$ can be constructed from the data already present!
Thus we make the following definitions:

\begin{defn} \label{coherenthomo} \et
A {\bf homomorphism} $F \maps C \to C'$ between coherent 2-groups
is a weak monoidal functor.
\end{defn}

\begin{defn} \label{coherent2homo} \et
A {\bf 2-homomorphism} $\theta \maps F \To G$ between
homomorphisms \break $F,G \maps C \to C'$ of coherent 2-groups
is a weak monoidal natural transformation.
\end{defn}

\begin{defn} \et Let {\bf C2G} be the 2-category consisting of coherent
2-groups, homomorphisms between these, and 2-homomorphisms between
those.  \end{defn}

\noindent It is clear that $\cg$ forms a 2-category since it
is actually a full and 2-full sub-2-category of $\Mon$.

Now let us describe how to define $F_{-1}$ in terms
of the other data in a coherent 2-group homomorphism $F \maps C \to C'$.
By analogy with $F_2$ and $F_0$, we would expect $F_{-1}$ to
make these diagrams commute:
\begin{description}
    \item[H1]
    \[
    \\
    \vcenter{
    \xymatrix@!C{
       F(x) \ten \overline{F(x)}
       \ar[r]^<<<<<{1 \ten F_{-1}}
     & F(x) \ten F(\xb)
          \ar[r]^{F_{2}}
     & F(x \ten \xb)   \\
        1'
          \ar[u]^{i_{F(x)}}
          \ar[rr]^{F_{0}}
     && F(1)
          \ar[u]_{F(i_x)}
     }}
\]
    \item[H2]
    \[
    \\
    \vcenter{
    \xymatrix@!C{
       \overline{F(x)} \ten F(x)
       \ar[r]^<<<<{F_{-1} \ten 1}
       \ar[d]_{e_{F(x)}}
     & F(\xb) \ten F(x)
          \ar[r]^{F_{2}}
     & F(\xb \ten x)
       \ar[d]^{F(e_x)} \\
        1'
          \ar[rr]^{F_{0}}
     && F(1)
     }}
\]
\end{description}
for all $x \in C$.  These diagrams say that $F_{-1}$ gets along with units
and counits.

Suppose we wish to construct an isomorphism
that simultaneously satisfies both of these coherence laws. To do this we
can take one of these coherence laws, solve it for $F_{-1}$, and prove
that the result automatically satisfies the other coherence law!
To do this, we start with the axiom \textbf{H1} expressed
in a more suggestive manner:
\[
\def\objectstyle{\scriptstyle}
\def\labelstyle{\scriptstyle}
\vcenter{
 \xymatrix@!C{
    &  F(x) \ten \overline{F(x)}
      \ar[rr]^{1 \ten F_{-1}}
      \ar[dl]_{i_{F(x)}^{-1}}
    && F(x) \ten F(\xb) \\
    1' \ar[drr]_{F_{0}}
    &&&& F(x \ten \xb) \ar[ul]_{F_{2}^{-1}} \\
    && F(1) \ar[urr]_{F(i_{x})}
    }}
\]
If we assume this diagram commutes, it gives a formula for
\[   1 \tensor F_{-1}
     \maps \xymatrix@1{F(x) \ten \overline{F(x)}
            \ar[r]^<<<<<{\sim} & \;     F(x) \ten F(\xb) .  }
\]

Next we shall solve for $F_{-1}$ by cancelling the tensor product $F(x)
\tensor $.  In the arguments that follow we will no longer include
subscripts on any of the structural maps $a, \ell, r, e$ or $i$ unless
confusion is likely to arise.  Tensoring the above morphism by
$\overline{F(x)}$ we obtain
\[   1 \tensor (1 \tensor F_{-1})
     \maps \xymatrix@1{
    \overline{F(x)} \ten (F(x) \ten \overline{F(x)})
            \ar[r]^<<<<<{\sim} & \;     F(x) \ten F(\xb)
    \overline{F(x)} \ten (F(x) \ten F(\xb))
  }
\]
However, the left-hand side is isomorphic to $\overline{F(x)}$ via
the following composite
\[
\def\objectstyle{\scriptstyle}
\def\labelstyle{\scriptstyle}
 \vcenter{
  \xy
    (-60,0)*{\overline{F(x)}}="1";
    (-40,0)*{1' \ten \overline{F(x)}}="2";
    (-2,0)*{(\overline{F(x)} \ten F(x)) \ten \overline{F(x)}}="3";
    (50,0)*{\overline{F(x)} \ten (F(x) \ten \overline{F(x)}) .}="4";
        {\ar^{\ell^{-1}}       "1";"2"};
        {\ar^<<<<<<<<<<{e^{-1} \ten 1}       "2";"3"};
        {\ar^{a}       "3";"4"};
 \endxy }
\]
Further, the right hand side is isomorphic to $F(\xb)$ via the
following composite:
\[
\def\objectstyle{\scriptstyle}
\def\labelstyle{\scriptstyle}
 \vcenter{
  \xy
    (50,0)*{F(\xb)  .}="1";
    (30,0)*{1' \ten F(\xb)}="2";
    (-10,0)*{(\overline{F(x)} \ten F(x)) \ten F(\xb)}="3";
    (-60,0)*{\overline{F(x)} \ten (F(x) \ten F(\xb))}="4";
        {\ar^{\ell}                          "2";"1"};
        {\ar^<<<<<<<<<<{e \ten 1}       "3";"2"};
        {\ar^{a} "4";"3"};
 \endxy }
\]
Stringing together these isomorphisms we obtain an isomorphism
from $\overline{F(x)}$ to $F(\xb)$ which is none other than
$F_{-1}$.  In other words, we have a commutative diagram
which serves to define $F_{-1}$:
\begin{description}
    \item[F1]
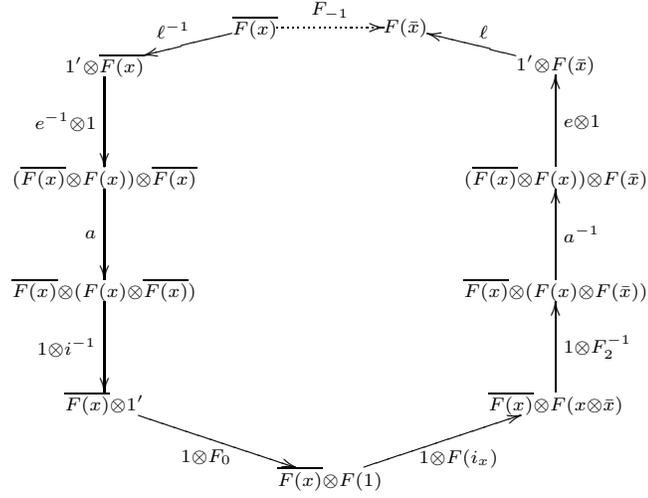
\[
\def\objectstyle{\scriptstyle}
\def\labelstyle{\scriptstyle}
 \vcenter{
 \xy
    (-10,30)*{\overline{F(x)}}="1";
    (-30,25)*{1' \ten \overline{F(x)}}="2";
    (-30,10)*{(\overline{F(x)} \ten F(x)) \ten \overline{F(x)}}="3";
    (-30,-5)*{\overline{F(x)} \ten (F(x) \ten \overline{F(x)})}="4";
    (-30,-20)*{\overline{F(x)} \ten 1'}="5";
    (0,-30)*{\overline{F(x)} \ten F(1)}="6";
    (30,-20)*{\overline{F(x)} \ten F(x \ten \xb)}="7";
    (30,-5)*{\overline{F(x)} \ten (F(x) \ten F(\xb))}="8";
    (30,10)*{(\overline{F(x)} \ten F(x)) \ten F(\xb)}="9";
    (30,25)*{1' \ten F(\xb)}="10";
    (10,30)*{F(\xb)}="11";
        {\ar_{\ell^{-1}} "1";"2"};
        {\ar_{e^{-1} \ten 1} "2";"3"};
        {\ar_{a} "3";"4"};
        {\ar_{1 \ten i^{-1}} "4";"5"};
        {\ar_{1 \ten F_{0}} "5";"6"};
        {\ar_{1 \ten F(i_x)} "6";"7"};
        {\ar_{1 \ten F_{2}^{-1}} "7";"8"};
        {\ar_{a^{-1}} "8";"9"};
        {\ar_{e \ten 1} "9";"10"};
        {\ar_{\ell} "10";"11"};
        {\ar@{.>}^{F_{-1}} "1";"11"};
 \endxy }
\]
\end{description}
We have derived this diagram \textbf{F1} from the assumption
\textbf{H1}; conversely one can show that \textbf{F1} implies
\textbf{H1}.

As a side remark, note that in the above diagram
could have mapped $\overline{F(x)}$ directly
to $\overline{F(x)} \ten 1'$ using $r^{-1}$.  This provides an
alternative definition of $F_{-1}$:

\begin{description}
    \item[${\bf F1'}$]
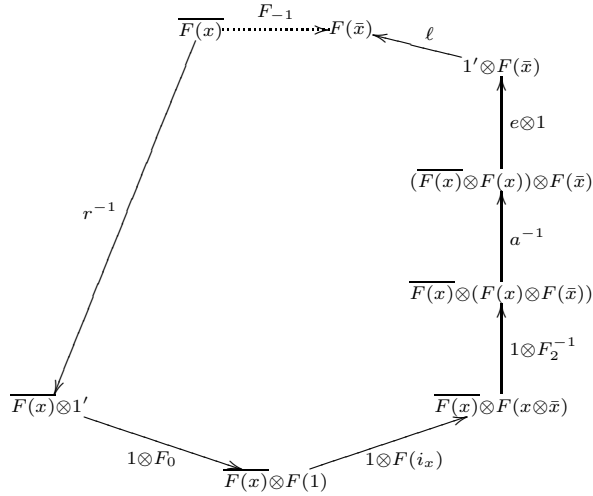
\[
\def\objectstyle{\scriptstyle}
\def\labelstyle{\scriptstyle}
 \vcenter{
 \xy
    (-10,30)*{\overline{F(x)}}="1";
    (-30,-20)*{\overline{F(x)} \ten 1'}="5";
    (0,-30)*{\overline{F(x)} \ten F(1)}="6";
    (30,-20)*{\overline{F(x)} \ten F(x \ten \xb)}="7";
    (30,-5)*{\overline{F(x)} \ten (F(x) \ten F(\xb))}="8";
    (30,10)*{(\overline{F(x)} \ten F(x)) \ten F(\xb)}="9";
    (30,25)*{1' \ten F(\xb)}="10";
    (10,30)*{F(\xb)}="11";
        {\ar_{r^{-1}}   "1";"5"};
        {\ar_{1 \ten F_{0}} "5";"6"};
        {\ar_{1 \ten F(i_x)} "6";"7"};
        {\ar_{1 \ten F_{2}^{-1}} "7";"8"};
        {\ar_{a^{-1}} "8";"9"};
        {\ar_{e \ten 1} "9";"10"};
        {\ar_{\ell} "10";"11"};
        {\ar@{.>}^{F_{-1}} "1";"11"};
 \endxy }
\]
\end{description}
The assertion that these definitions agree is equivalent to the
fact that this diagram commutes:
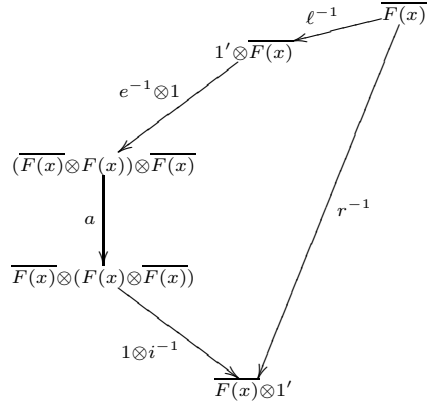
\[
\def\objectstyle{\scriptstyle}
\def\labelstyle{\scriptstyle}
 \vcenter{
 \xy
    (-10,30)*{\overline{F(x)}}="1";
    (-30,25)*{1' \ten \overline{F(x)}}="2";
    (-50,10)*{(\overline{F(x)} \ten F(x)) \ten \overline{F(x)}}="3";
    (-50,-5)*{\overline{F(x)} \ten (F(x) \ten \overline{F(x)})}="4";
    (-30,-20)*{\overline{F(x)} \ten 1'}="5";
        {\ar_{\ell^{-1}} "1";"2"};
        {\ar_{e^{-1} \ten 1} "2";"3"};
        {\ar_{a} "3";"4"};
        {\ar_{1 \ten i^{-1}} "4";"5"};
        {\ar^{r^{-1}}    "1";"5"};
 \endxy }
\]
Careful inspection reveals that this diagram is
none other than the second zig-zag identity satisfied by the
adjoint equivalence $(x,\xb,i_x,e_x)$!  So, there is no problem
of conflicting choices here.  In fact, this is guaranteed by
Ulbrich and Laplaza's coherence theorem for categories with group
structure \cite{Laplaza,Ulbrich}.

Alternatively, we could have defined $F_{-1}$ so that
\textbf{H2} is satisfied.  Then this diagram commutes:
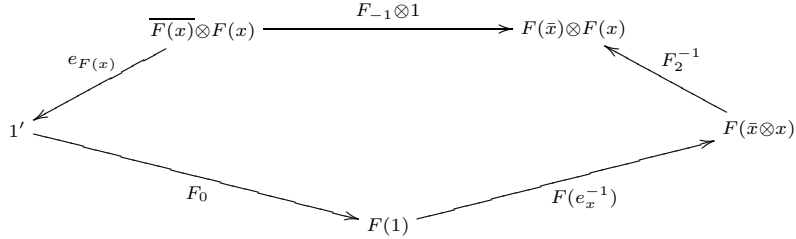
\[
\def\objectstyle{\scriptstyle}
\def\labelstyle{\scriptstyle}
\vcenter{
 \xymatrix@!C{
    &  \overline{F(x)} \ten F(x)
      \ar[rr]^{F_{-1} \ten 1}
      \ar[dl]_{e_{F(x)}}
    && F(\xb) \ten F(x)  \\
    1' \ar[drr]_{F_{0}}
    &&&& F(\xb \ten x) \ar[ul]_{F_{2}^{-1}} \\
    && F(1) \ar[urr]_{F(e_{x}^{-1})}
    }}
\]
Tensoring on the right by $\overline{F(x)}$ and applying a similar
argument as before we obtain two
additional possibilities for $F_{-1}$, which again are
actually equal thanks to the second zig-zag identity.  We denote
this way of defining $F_{-1}$ simply as \textbf{F2}:

\begin{description}
    \item[F2]
\[
\def\objectstyle{\scriptstyle}
\def\labelstyle{\scriptstyle}
 \vcenter{
 \xy
    (-10,30)*{\overline{F(x)}}="1";
    (-30,25)*{ \overline{F(x)} \ten1'}="2";
    (-30,10)*{\overline{F(x)} \ten (F(x) \ten \overline{F(x)})}="3";
    (-30,-5)*{(\overline{F(x)} \ten F(x)) \ten \overline{F(x)}}="4";
    (-30,-20)*{1' \ten \overline{F(x)}}="5";
    (0,-30)*{F(1) \ten \overline{F(x)}}="6";
    (30,-20)*{F(\xb \ten x) \ten \overline{F(x)}}="7";
    (30,-5)*{(F(\xb) \ten F(x)) \ten \overline{F(x)}}="8";
    (30,10)*{F(\xb) \ten (F(x) \ten \overline{F(x)})}="9";
    (30,25)*{F(\xb) \ten 1'}="10";
    (10,30)*{F(\xb)}="11";
        {\ar_{r^{-1}}        "1";"2"};
        {\ar_{ 1 \ten i}  "2";"3"};
        {\ar_{a^{-1}}                          "3";"4"};
        {\ar_{e \ten  1}    "4";"5"};
        {\ar_{F_{0} \ten 1}      "5";"6"};
        {\ar_{F(e_x^{-1}) \ten 1}"6";"7"};
        {\ar_{F_{2}^{-1} \ten  1}"7";"8"};
        {\ar_{a}                               "8";"9"};
        {\ar_{1 \ten i}^{-1}     "9";"10"};
        {\ar_{r}                      "10";"11"};
        {\ar@{.>}^{\ell^{-1}} "1";"5"};
        {\ar@{.>}^{F_{-1}} "1";"11"};
 \endxy }
\]
\end{description}

In fact, just as \textbf{F1} is equivalent to \textbf{H1}, \textbf{F2}
is equivalent to \textbf{H2}.  Next, we would like to show that
\textbf{F1} and \textbf{F2} give the same definition of $F_{-1}$, so
that using either to define this isomorphism guarantees that both
\textbf{H1} and \textbf{H2} are satisfied.  To accomplish this we shall
assume \textbf{H1} and use this to establish \textbf{F2}.  Consider the
following diagram:
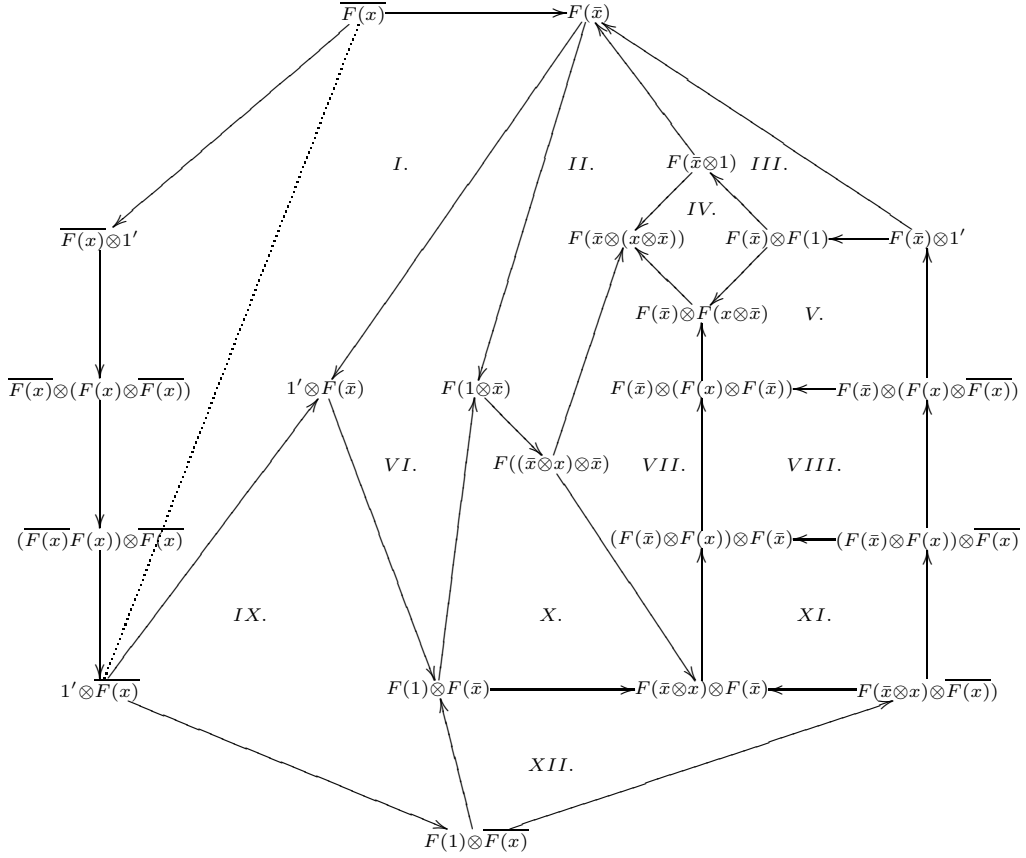
\[
\def\objectstyle{\scriptstyle}
\def\labelstyle{\scriptstyle}
 \vcenter {
 \xy
    (-25,60)*{\overline{F(x)}}="1";
    (-60,30)*{\overline{F(x)} \ten 1'}="2";
    (-60,10)*{\overline{F(x)} \ten (F(x) \ten \overline{F(x)})}="3";
    (-60,-10)*{(\overline{F(x)} F(x)) \ten \overline{F(x)}}="4";
    (-60,-30)*{1' \ten \overline{F(x)}}="5";
    (-10,-50)*{F(1) \ten \overline{F(x)}}="6";
    (50,-30)*{F(\xb \ten x) \ten \overline{F(x)})}="7";
    (50,-10)*{(F(\xb) \ten F(x)) \ten \overline{F(x)}}="8";
    (50,10)*{F(\xb) \ten (F(x) \ten \overline{F(x)})}="9";
    (50,30)*{F(\xb) \ten 1'}="10";
    (5,60)*{F(\xb)}="11";
            (-15,-30)*{F(1) \ten F(\xb)}="J";
            (20,-30)*{F(\xb \ten x) \ten F(\xb)}="7'";
            (20,-10)*{(F(\xb) \ten F(x)) \ten F(\xb)}="H";
            (20,10)*{F(\xb) \ten (F(x) \ten F(\xb))}="G";
            (20,20)*{F(\xb) \ten F(x \ten \xb)}="A";
            (30,30)*{F(\xb) \ten F(1)}="B";
            (20,40)*{F(\xb \ten 1)}="C";
            (10,30)*{F(\xb \ten (x \ten \xb))}="D";
            (-10,10)*{F(1 \ten \xb)}="F";
            (0,0)*{F((\xb \ten x) \ten \xb)}="E";
            (-30,10)*{1' \ten F(\xb)}="10'";
        {\ar^{}       "1";"11"};
        {\ar^{}       "1";"2"};
        {\ar^{}       "2";"3"};
        {\ar^{}       "3";"4"};
        {\ar^{}       "4";"5"};
        {\ar^{}       "5";"6"};
        {\ar^{}       "6";"7"};
        {\ar^{}       "7";"8"};
        {\ar^{}       "8";"9"};
        {\ar^{}       "9";"10"};
        {\ar^{}       "10";"11"};
            {\ar@{.}^{}       "1";"5"};
            {\ar^{}       "7";"7'"};
            {\ar^{}       "J";"7'"};
            {\ar^{}       "7'";"H"};
            {\ar^{}       "6";"J"};
            {\ar^{}       "8";"H"};
            {\ar^{}       "H";"G"};
            {\ar^{}       "9";"G"};
            {\ar^{}       "G";"A"};
            {\ar^{}       "B";"A"};
            {\ar^{}       "10";"B"};
            {\ar^{}       "B";"C"};
            {\ar^{}       "C";"D"};
            {\ar^{}       "A";"D"};
            {\ar^{}       "C";"11"};
            {\ar^{}       "11";"F"};
            {\ar^{}       "F";"E"};
            {\ar^{}       "E";"D"};
            {\ar^{}       "11";"10'"};
            {\ar^{}       "5";"10'"};
            {\ar^{}       "10'";"J"};
            {\ar^{}       "J";"F"};
            {\ar^{}       "E";"7'"};
        (-20,40)*{I.};
        (4,40)*{II.};
        (29,40)*{III.};
        (20,34)*{IV.};
        (35,20)*{V.};
        (-20,0)*{VI.};
        (15,0)*{VII.};
        (35,0)*{VIII.};
        (-40,-20)*{IX.};
        (0,-20)*{X.};
        (35,-20)*{XI.};
        (0,-40)*{XII.};
     \endxy
 }
\]
Squares $I$, $IV$, $VIII$, $X$, $XI$, $XI$, and $XII$ commute from the
naturality of the isomorphisms $\ell$, $F_{2}$, $a$, $e_x$. Application of
$F$ to the second zig-zag law gives $II$.  Diagrams $III$, $VI$ and
$VII$ commute by the definition of weak monoidal functor, while $IX$
commutes by a well-known property of monoidal categories.
Diagram $V$ is merely \textbf{H1} tensored on the left
by $F(\xb)$. Thus, our choice of $F_{-1}$ satisfies both
\textbf{H1} and \textbf{H2}.  It follows that $F_{-1}$ and its
coherence laws are superfluous to the definition of a coherent
2-group homomorphism; we get them `for free'.

\section{Internalization} \label{internalizationsection}

`Internalization' is valuable tool for generalizing concepts from
the category of sets to other categories.  To internalize a
concept, we need to express it in a purely diagrammatic form. As
an example, consider the ordinary notion of group.  We can define
this notion using commutative diagrams by specifying:
\begin{description}
    \item[(i)] a set $G$,
    \item[(ii)]  a multiplication function $m \maps G \times G \to G$,
    \item[(iii)] a unit for the multiplication given by the function
    $e \maps 1 \to G$ where $1$ is the terminal object in $\Set$,
    \item[(iv)] a function $\inv \maps G \to G$,
\end{description}
making the following diagrams commute:
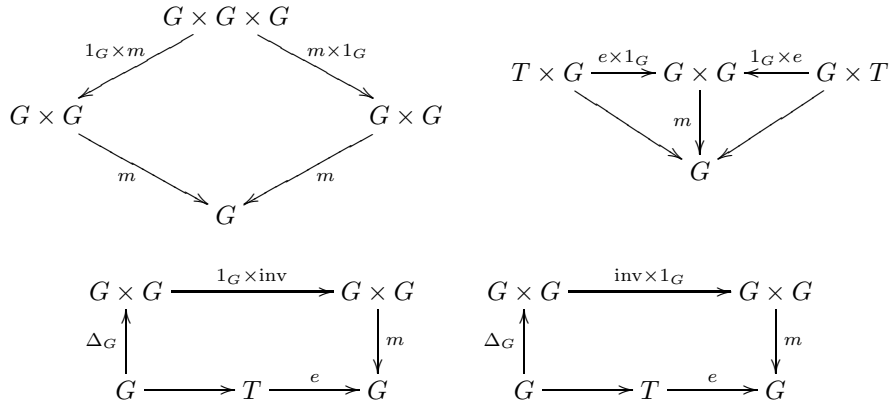
\[ \vcenter{
  \xymatrix{
&   G \times G \times G \ar[dr]^{m \times 1_G}
      \ar[dl]_{1_G \times m} \\
    G \times G \ar[dr]_{m}
&&  G \times G \ar[dl]^{m}  \\
&  G
  }}
\qquad \vcenter{
  \xymatrix{
    T \times G \ar[r]^{e \times 1_G} \ar[dr]
  & G \times G \ar[d]_{m}
  & G \times T \ar[l]_{1_G \times e} \ar[dl] \\
  & G
  }}
\]
\[
  \xymatrix{
     G \times G \ar[rr]^{1_G \times \inv}
  && G \times G \ar[d]^{m} \\
     G \ar[u]^{\Delta_G} \ar[r]
  &  T \ar[r]^{e}
  &  G
  }
  \qquad
  \xymatrix{
     G \times G \ar[rr]^{ \inv \times 1_G}
  && G \times G \ar[d]^{m} \\
     G \ar[u]^{\Delta_G} \ar[r]
  &  T \ar[r]^{e}
  &  G
  }
\]
where $\Delta_G$ is the diagonal map.  To internalize this concept one
replaces the set $G$ by an object in an arbitrary category $C$ with
finite products, and the functions $m, e,$ and $inv$ by morphisms in
$C$.  Making these substitutions in the definition above we arrive at
the definition of a {\bf group object in $C$}.  An ordinary group is the
special case where $C = \Set$.  We can also define a strict 2-group to
be a group object in $\Cat$.  Similarly, a topological group is
a group object in $\Top$, and a Lie group is a group object in $\Diff$.

We would like to define a `coherent 2-group object' in a similar manner.
To motivate this definition it is helpful to use `string diagrams'
\cite{JS0,Street}.  These are Poincar\'e dual to the globular diagrams
previously used in the theory of bicategories.  In other words, to
obtain a string diagram one draws objects as 2-dimensional regions
in the plane, 1-morphisms as 1-dimensional `strings' separating regions,
and 2-morphisms as 0-dimensional points (or small balls, if we wish to
label them).  We will only need these diagrams in the special case of a weak
monoidal category, which we will think of as a bicategory with a single
object, say $\bullet$.  A morphism $f \maps x \to y $ in a weak monoidal
category corresponds to a 2-morphism in a bicategory with one object,
and we convert the globular picture of this into a string diagram as
follows:
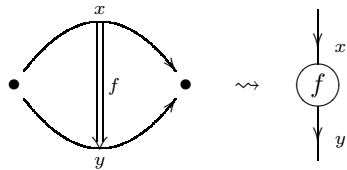
\[
 \xymatrix{
    \bullet
      \ar@/^2pc/[rr]_{\quad}^{x}="1"
      \ar@/_2pc/[rr]_{y}="2"
  && \bullet
    \ar@{=>}"1" ;"2"^{f}
  }
\quad \rightsquigarrow \quad
 \xy
 (0,0)*++{f}*\cir{}="f";
 (0,10)**\dir{-} ?(.5)*\dir{<}+(3,0)*{\scriptstyle x};
 "f";(0,-10)**\dir{-} ?(.75)*\dir{>}+(3,0)*{\scriptstyle y};
\endxy
\]
where diagrams are read from top to bottom.  Composition of
morphisms is achieved by placing one of the strings on top of the
other. For instance, given $f \maps x \to y$ and $g \maps y \to z$,
their composite is depicted as the following:
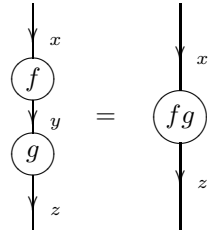
\[
\xy
 (0,5)*++{f}*\cir{}="f";
 (0,15)**\dir{-} ?(.5)*\dir{<}+(3,0)*{\scriptstyle x};
 "f";(0,-5)*++{g}*\cir{}="g";
 **\dir{-} ?(.4)*\dir{<}+(3,0)*{\scriptstyle y};
 "f";"g";(0,-15)**\dir{-} ?(.75)*\dir{>}+(3,0)*{\scriptstyle z};
\endxy
\quad = \quad
\xy
 (0,15)*{};
 (0,0)*++{fg}*\cir{}="fg";
 **\dir{-} ?(.5)*\dir{<}+(3,0)*{\scriptstyle x};
 "fg";(0,-15)*{};
 **\dir{-} ?(.4)*\dir{<}+(3,0)*{\scriptstyle z};
\endxy
\]
This diagram is Poincar\'e dual to the globular way of drawing
composition of 2-morphisms in a bicategory:
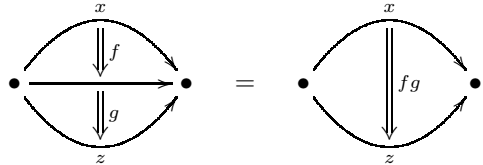
\[
 \xymatrix{
    \bullet
      \ar@/^2pc/[rr]^{x}_{}="0"
      \ar[rr]^{}="1"
      \ar@/_2pc/[rr]_{z}_{}="2"
      \ar@{=>}"0";"1" ^{f}
      \ar@{=>}"1";"2" ^{g}
 &&  \bullet
 }
 \quad = \quad
  \xymatrix{
    \bullet
      \ar@/^2pc/[rr]^{x}_{}="0"
      \ar@/_2pc/[rr]_{z}_{}="2"
      \ar@{=>}"0";"2" ^{fg}
 &&  \bullet
 }
\]

Tensoring objects in the monoidal category will be written by
setting arrows side by side; the unit object will not be drawn
in the diagrams, but merely implied. As an example of this,
consider how we obtain the string diagram
corresponding to $i_{x} \maps 1 \to x \ten \xb$:
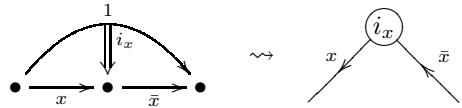
\[
\vcenter{
 \xymatrix{
    \bullet \ar[r]_{x} \ar@/^2pc/[rr]^{1}="1"
    & \bullet \ar[r]_{\xb} \ar@{=>}^<<<{i_{x}}"1";{}
    & \bullet
 } }
\quad \rightsquigarrow \quad
 \vcenter {\xy
    (0,0)*+{i_{x}} *\cir{}="i";
    (-10,-10)*{};
    **\dir{-} ?(.4)*\dir{<}+(-1,2)*{\scriptstyle x};
    (10,-10)*{};"i";
    **\dir{-} ?(.5)*\dir{<}+(3,1)*{\scriptstyle \xb};
 \endxy }
\]
Note that weak inverse objects are written as arrows `going
backwards in time'.  We will find it beneficial to draw the above
diagram in a simpler form:
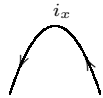
\[
\xy
    (-6,0)*{};
    (6,0)*{};
      **\crv{(0,18)} ?(.16)*\dir{>} ?(.9)*\dir{>};
    (1,11)*{\scriptstyle i_{x}};
\endxy
\]
where it is understood that the downward pointing arrow
corresponds to $x$ and the upward pointing arrow to
$\xb$.  Similarly, we draw the morphism $e_x$ as
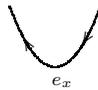
\[
 \xy
    (-6,4)*{};(6,4)*{};
      **\crv{(0,-12)}
       ?(.20)*\dir{>}+(2,-1)
       ?(.85)*\dir{>}+(-2,-1);
      (1,-6)*{\scriptstyle e_{x}};
\endxy
\]
In this notation, the zig-zag identities become

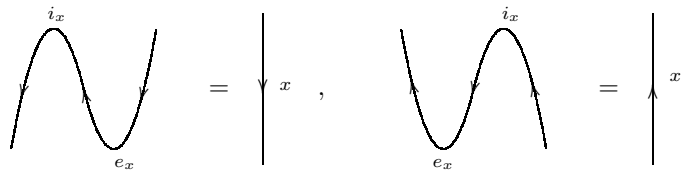
\[
\xy
    (0,-10)*++{}="g";
    (10,0)*{}="mid";
      **\crv{(5,20)}  ?(.83)*\dir{>};
    (6.5,10)*{\scriptstyle i_{x}};
    (20,10)*++{}="f";
     "f";"mid"; **\crv{(15,-20)} ?(0)*\dir{<} ?(.76)*\dir{<};
    (15.5,-10)*{\scriptstyle e_{x}};
\endxy
\quad = \quad
 \xy
  (0,10)*{};(0,-10)*{};
  **\dir{-} ?(.47)*\dir{<}+(3,1)*{\scriptstyle x}
 \endxy
\quad , \qquad
 \xy
    (0,10)*++{}="g";
    (10,0)*{}="mid";
      **\crv{(5,-20)} ?(.02)*\dir{>} ?(.83)*\dir{>};
    (6,-10)*{\scriptstyle e_{x}};
    (20,-10)*++{}="f";
     "f";"mid"; **\crv{(15,20)}  ?(.78)*\dir{<};
    (15,10)*{\scriptstyle i_{x}};
\endxy
\quad = \quad
 \xy
  (0,10)*{};(0,-10)*{};
  **\dir{-} ?(.53)*\dir{>}+(3,1)*{\scriptstyle x}
 \endxy
\]
which explains their name.

We would like to define a coherent 2-group using only commutative
diagrams, so that the groundwork will have been laid for the
definition of the more general notion of  a `coherent group
object' in a 2-category with products.  To begin, notice that
for any coherent 2-group $C$ there is a functor
${}^{-1} \maps C^{\rm op} \to C$, expressed diagrammatically as
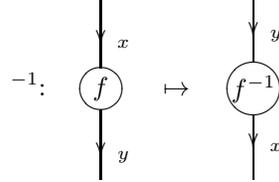
\[ ^{-1} \maps \quad
\xy
 (0,0)*++{f}*\cir{}="f";
 (0,12)**\dir{-} ?(.5)*\dir{<}+(3,0)*{\scriptstyle x};
 "f";(0,-12)**\dir{-} ?(.75)*\dir{>}+(3,0)*{\scriptstyle y};
\endxy
\quad \mapsto \quad
 \xy
 (0,0)*+{f^{-1}}*\cir{}="f";
 (0,12)**\dir{-} ?(.59)*\dir{<}+(3,0)*{\scriptstyle y};
 "f";(0,-12)**\dir{-} ?(.65)*\dir{>}+(3,0)*{\scriptstyle x};
\endxy
\]
and coming from the invertibility of morphisms in a coherent
2-group.    There is also a functor $^{*} \maps C^{\rm op} \to C$ sending
each object $x \in C$ to its specified weak inverse $\bar{x}$,
and acting on morphisms as follows:
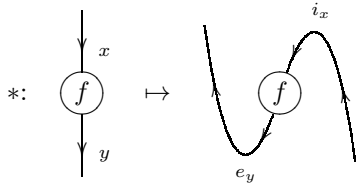
\[ * \maps \quad
\xy
 (0,0)*++{f}*\cir{}="f";
 (0,11)**\dir{-} ?(.5)*\dir{<}+(3,0)*{\scriptstyle x};
 "f";(0,-11)**\dir{-} ?(.75)*\dir{>}+(3,0)*{\scriptstyle y};
\endxy
\quad  \mapsto \quad
 \xy
    (0,0)*++{f}*\cir{}="f";
    (10,-9)*{};
      **\crv{(6,20)} ?(.2)*\dir{>} ?(.86)*\dir{>};
      "f";(-10,9)*{}; **\crv{(-6,-20)} ?(.15)*\dir{<} ?(.78)*\dir{<};
    (5.5,11)*{\scriptstyle i_{x}};
    (-4.5,-11)*{\scriptstyle e_{y}}
 \endxy
\]

Composing these two contravariant functors $^{-1}$ and $*$ we
construct a covariant functor $\inv \maps C \to C$ given by
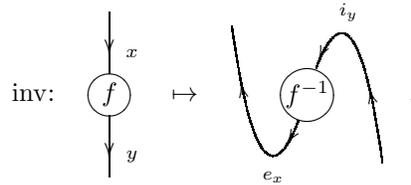
\[ \inv \maps \quad
\xy
 (0,0)*++{f}*\cir{}="f";
 (0,11)**\dir{-} ?(.5)*\dir{<}+(3,0)*{\scriptstyle x};
 "f";(0,-11)**\dir{-} ?(.75)*\dir{>}+(3,0)*{\scriptstyle y};
\endxy
\quad  \mapsto \quad
 \xy
    (0,0)*+{f^{-1}}*\cir{}="f";
    (10,-9)*{};
      **\crv{(6,20)} ?(.2)*\dir{>} ?(.86)*\dir{>};
      "f";(-10,9)*{}; **\crv{(-6,-20)} ?(.15)*\dir{<} ?(.78)*\dir{<};
    (5.5,11)*{\scriptstyle i_{y}};
    (-4.5,-11)*{\scriptstyle e_{x}}
 \endxy
 \quad .
\]
In order to prove the functoriality of $\inv$,
consider two composable morphisms $ f \maps x \to y$ and $g \maps
y \to z$.  The equation $ \inv(fg) = \inv(f) \inv(g)$ becomes
the following in string diagram notation:
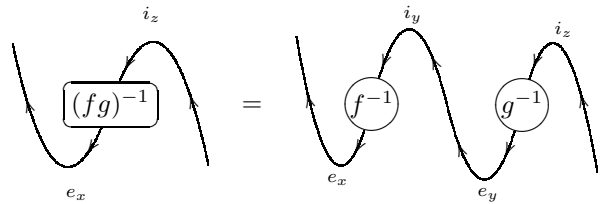
\[
 \xy
    (0,0)*+{(fg)^{-1}}*\frm<3pt>{-}="f";
    (13,-8)*{};
      **\crv{(7,20)} ?(.2)*\dir{>} ?(.86)*\dir{>};
      "f";(-13,8)*{}; **\crv{(-7,-20)} ?(.15)*\dir{<} ?(.78)*\dir{<};
    (5.5,12)*{\scriptstyle i_{z}};
    (-4.5,-12)*{\scriptstyle e_{x}}
 \endxy
\quad = \quad
 \xy
    (0,0)*+{f^{-1}}*\cir{}="g";
    (10,0)*{}="mid";
      **\crv{(5,20)} ?(.2)*\dir{>} ?(.86)*\dir{>};
      "g";(-10,9)*{}; **\crv{(-5,-20)} ?(.15)*\dir{<} ?(.78)*\dir{<};
    (5.5,12)*{\scriptstyle i_{y}};
    (-4.5,-10)*{\scriptstyle e_{x}};
    (20,0)*+{g^{-1}}*\cir{}="f";
     "f";"mid"; **\crv{(15,-20)} ?(.15)*\dir{<} ?(.78)*\dir{<};
     "f";(30,-9) **\crv{(25,20)} ?(.15)*\dir{<} ?(.78)*\dir{<};
    (15.5,-12)*{\scriptstyle e_{y}};
    (25.5,10)*{\scriptstyle i_{z}};
 \endxy
\]
Thus, in order for equality to hold we must have
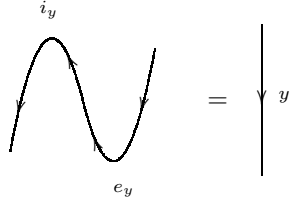
\[
\xy
    (0,-9)*++{}="g";
    (10,0)*{}="mid";
      **\crv{(5,20)} ?(.2)*\dir{>} ?(.86)*\dir{>};
    (5.5,12)*{\scriptstyle i_{y}};
    (20,9)*++{}="f";
     "f";"mid"; **\crv{(15,-20)} ?(.15)*\dir{<} ?(.78)*\dir{<};
    (15.5,-12)*{\scriptstyle e_{y}};
\endxy
\quad = \quad
 \xy
  (0,10)*{};(0,-10)*{};
  **\dir{-} ?(.47)*\dir{<}+(3,1)*{\scriptstyle y}
 \endxy \quad
\]
This diagram is merely the second zig-zag identity!

Now, when all the unit and counit are isomorphisms either of the zig-zag
identities implies the other \cite{Rivano}.  In this case,
assuming that $\inv$ is a functor implies {\it both} zig-zag
identities.  This observation allows us to give an
equivalent definition of `coherent 2-group':

\begin{defn} \label{defc}  \et
A { \bf coherent 2-group} consists of:
 \begin{description}
    \item[(i)] a category $C$,
    \item[(ii)] a functor $m \maps C\times C \to C$, where we
    write $m(x,y)=x \ten y$ and $m(f,g)=f \ten g$ for objects $x, y,
    \in C$ and morphisms $f, g$ in $C$,
    \item[(iii)] a functor $\id \maps I \to C$ where $I$ is the
    terminal category, and we write the object in the range of
    this functor as $1 \in C$,
    \item[(iv)] a functor $ \inv \maps C \to C$,
    \item[(v)] natural isomorphisms $a,\ell,r,e,i$ as follows:
    \item[]
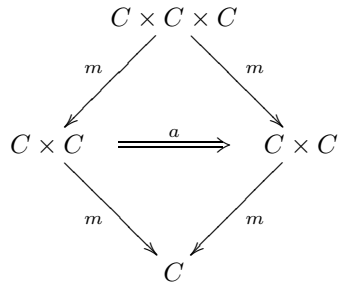
\[
\xymatrix@!0{
    && C \times C \times C
        \ar[ddrr]^m \ar[ddll]_m \\
            \\
    C \times C
        \ar[ddrr]_m
     & \ar@{=>}[rr]^{a}
     &&
     & C \times C
        \ar[ddll]^m \\
            \\
    && C  }
\]
    \item[]
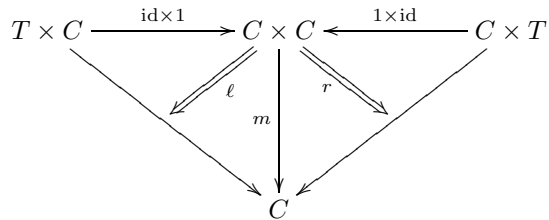
\[
 \xymatrix{
    T \times C
        \ar[rr]^{\id \times 1}
        \ar[ddrr]
     && C \times C
        \ar[dd]_m
        \ar@{=>}[dr]_{r}
        \ar@{=>}[dl]^{\ell}
     && C \times T
        \ar[ll]_{1 \times \id}
        \ar[ddll] \\
    & && \\
    && C  }
\]
    \item[]
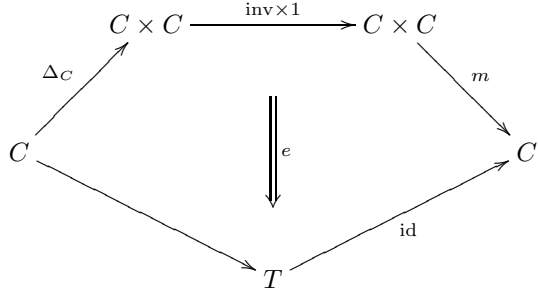
\[
 \xymatrix@!0{
     && C \times C
        \ar[rrrr]^{\inv \times 1}
      &&&& C \times C
        \ar[ddrr]^{m}  \\
     &&&&
        \ar@{=>}[dd]^e \\
    C
        \ar[urur]^{\Delta_C}
        \ar[ddrrrr]
     &&&&&&&& C \\
    &&&&   \\
    &&&& T \ar[uurrrr]_{\id}    }
\]

    \item[]
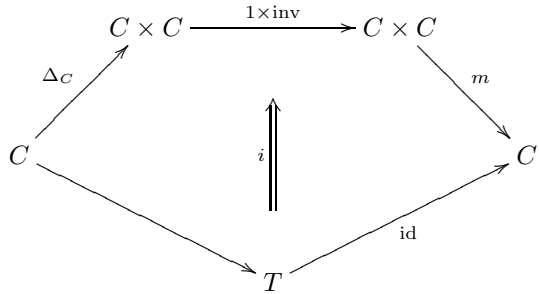
\[
 \xymatrix@!0{
     && C \times C
        \ar[rrrr]^{ 1 \times \inv}
      &&&& C \times C
        \ar[ddrr]^{m}  \\
     &&&&  \\
    C
        \ar[urur]^{\Delta_C}
        \ar[ddrrrr]
     &&&&&&&& C \\
    &&&& \ar@{=>}[uu]^i  \\
    &&&& T \ar[uurrrr]_{\id}    }
\]
making the following diagrams commute:
    \item[A1]
\[
 \xy
    (0,20)*{(x \ten y)(z \ten w)}="1";
    (40,0)*{x \ten (y \ten (z \ten w))}="2";
    (25,-20)*{ \quad x \ten ((y \ten z) \ten w)}="3";
    (-25,-20)*{(x \ten (y \ten z)) \ten w}="4";
    (-40,0)*{((x \ten y) \ten z) \ten w}="5";
        {\ar^{a_{x,y,z \ten w}}     "1";"2"}
        {\ar_{1_x \ten a _{y,z,w}}  "3";"2"}
        {\ar^{a _{x,y \ten z,w}}    "4";"3"}
        {\ar_{a _{x,y,z} \ten 1_w}  "5";"4"}
        {\ar^{a _{x \ten y,z,w}}    "5";"1"}
 \endxy
 \\
\]

    \item[A2]
\[
 \xymatrix{
    (x \ten 1) \ten y
        \ar[rr]^{a_{x,1,y}}
        \ar[dr]_{r_x \ten y}
     &&  x \ten (1 \ten y)
        \ar[dl]^{x \ten \ell_y } \\
    & x \ten y   } \\
\]

    \item[A3]
\[
 \xymatrix@!C{
     1 \ten x \ar[r]^<<<<<<{i_x \ten 1} \ar[d]_{\ell_x}
     & (x \ten \xb) \ten x \ar[r]^{a_{x, \xb ,x}}
     & x \ten ( \xb \ten x) \ar[d]^{1 \ten e_x} \\
     x \ar[rr]_{r^{-1}_x}
      && x \ten 1  } \\
\]

    \item[A4]
\[
 \xymatrix@!C{
    \xb \ten 1
        \ar[r]^<<<<<<{1 \ten i_x}
        \ar[d]_{r_{\xb}}
      & \xb \ten (x \ten \xb)
        \ar[r]^{a^{-1}_{\xb, x, \xb}}
      & (\xb \ten x) \ten \xb
        \ar[d]^{e_x \ten 1} \\
    \xb
        \ar[rr]_{\ell^{-1}_{\xb}}
      && 1 \ten \xb   } \\
\]
\end{description}
\end{defn}

We are now almost ready to define a `coherent 2-group object' in an
arbitrary 2-category with finite products.  Before we can do this, we
must be careful to phrase the definition in a way that does not use any
properties of a particular 2-category and its objects (in this case
$\Cat$).  Rather, we must formulate the definition in a purely
2-categorical way.  We have not done this yet in Definition \ref{defc},
since we made explicit use of the fact that $C$ was a category: the
commutative diagrams for the natural isomorphisms have \textit{objects}
of $C$ labelling their vertices.  We must correct this oversight.  The
problem, well-known to experts, is that making explicit mention of the
objects of $C$ has suppressed one dimension in the coherence diagrams.
This can be seen in axiom {\bf A5}:
\[
 \xy
    (0,20)*{(x \ten y)(z \ten w)}="1";
    (40,0)*{x \ten (y \ten (z \ten w))}="2";
    (25,-20)*{ \quad x \ten ((y \ten z) \ten w) .}="3";
    (-25,-20)*{(x \ten (y \ten z)) \ten w}="4";
    (-40,0)*{((x \ten y) \ten z) \ten w}="5";
        {\ar^{a_{x,y,z \ten w}}     "1";"2"}
        {\ar_{1_x \ten a _{y,z,w}}  "3";"2"}
        {\ar^{a _{x,y \ten z,w}}    "4";"3"}
        {\ar_{a _{x,y,z} \ten 1_w}  "5";"4"}
        {\ar^{a _{x \ten y,z,w}}    "5";"1"}
 \endxy
\]
The 2-dimensional appearance of this diagram results from
mentioning the objects $x,y,z,w \in C$.
We can avoid this by working with (for example) the
functor $(1 \times 1 \times m) \circ (1 \times m) \circ m$
instead of its value on the object $(x,y,z,w) \in C^4$, namely
$x \tensor (y \tensor (z \tensor w))$.  If we do this,
we see that the diagram is actually 3-dimensional!  It is
a pentagonal prism, a bit difficult to draw:
\\
\[
\xymatrix@!0{
    &&& C \times C \times C \times C \\
    &   \ar[ddddd]
    &   \ar[dddddd]
    &   \ar@{-}[ddd]
    &   \ar[dddddd]
    &   \ar[ddddd]           \\
    &&& \ar@{=>}[drr]|\hole  \\
    &   \ar@{=>}[urr]|\hole
        \ar@{=>}[dr]
    &&&&                     \\
    &&  \ar@{=>}[rr]
    &   \ar[d]
    &   \ar@{=>}[ur]         \\
    &&&                      \\
    & &&&&                   \\
    && & C &   }
\\
\]
where the downwards-pointing single arrows
are functors from $C^4$ to $C$, and the horizontal double arrows
are natural transformations between these functors, forming
a commutative pentagon.  Luckily we can also draw this pentagon
in a 2-dimensional way again, as follows:
\[
\def\objectstyle{\scriptstyle}
\def\labelstyle{\scriptstyle}
 \xy
    (-45,0)*{(m \ten 1 \ten 1) \circ (m \ten 1) \circ m}="1";
    (-45,-13)*{ (1 \ten m \ten 1)  \circ (m \ten 1) \circ m}="2";
    (45,-13)*{ (1 \ten m \ten 1) \circ ( 1 \ten m ) \circ m}="3";
    (45,0)*{(1 \ten 1 \ten m) \circ (1 \ten m) \circ m }="4";
    (0,0)*{(m \ten m) \circ m}="5";
    (-8,0)*{}="5'";
        {\ar^{a \ten 1 } "1";"2"};
        {\ar^>>>>>>>>>>>>>{(m \ten 1 \ten 1) \circ a} "1";"5"};
        {\ar_{(1 \ten m \ten 1) \circ a} "2";"3"};
        {\ar^{1 \ten a } "3";"4"};
        {\ar^<<<<<<<<<<<<<{(1 \ten 1 \ten m) \circ a} "5";"4"};
 \endxy
\]
This style of writing coherence laws allows us to make the
following definition:

\begin{defn} \label{coherentobject} \et
A { \bf coherent 2-group object} in a 2-category $K$ with finite
products consists of
\begin{description}
    \item[(i)] an object $C \in K$,
    \item[(ii)] a morphism $m \maps C\times C \to C$,
    \item[(iii)] a morphism $\id \maps I \to C$ where $I$ is the
    terminal object of $K$,
    \item[(iv)] a morphism $ \inv \maps C \to C$,
    \item[(v)] 2-isomorphisms $a,\ell,r,e,i$ as follows:
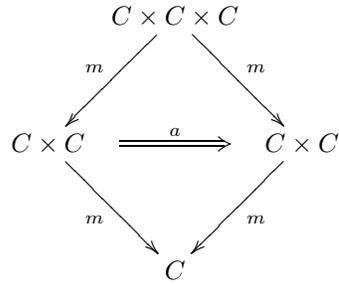
\[
\xymatrix@!0{
    && C \times C \times C
        \ar[ddrr]^m \ar[ddll]_m \\
            \\
    C \times C
        \ar[ddrr]_m
     & \ar@{=>}[rr]^{a}
     &&
     & C \times C
        \ar[ddll]^m \\
            \\
    && C  }
\]
    \item[]
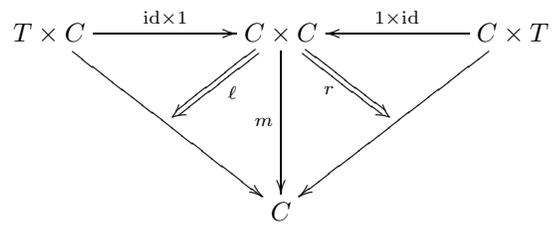
\[
 \xymatrix{
    T \times C
        \ar[rr]^{\id \times 1}
        \ar[ddrr]
     && C \times C
        \ar[dd]_m
        \ar@{=>}[dr]_{r}
        \ar@{=>}[dl]^{\ell}
     && C \times T
        \ar[ll]_{1 \times \id}
        \ar[ddll] \\
    & && \\
    && C  }
\]
    \item[]
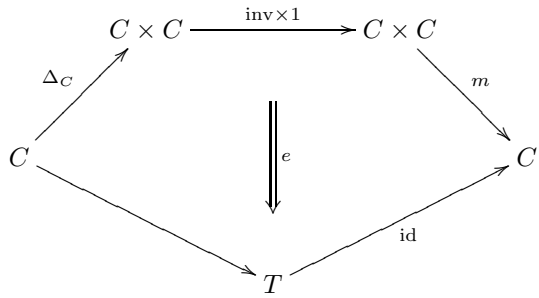
\[
 \xymatrix@!0{
     && C \times C
        \ar[rrrr]^{\inv \times 1}
      &&&& C \times C
        \ar[ddrr]^{m}  \\
     &&&&
        \ar@{=>}[dd]^e \\
    C
        \ar[urur]^{\Delta_C}
        \ar[ddrrrr]
     &&&&&&&& C \\
    &&&&   \\
    &&&& T \ar[uurrrr]_{\id}    }
\]

    \item[]
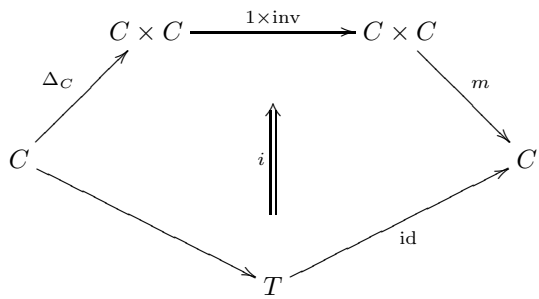
\[
 \xymatrix@!0{
     && C \times C
        \ar[rrrr]^{ 1 \times \inv}
      &&&& C \times C
        \ar[ddrr]^{m}  \\
     &&&&  \\
    C
        \ar[urur]^{\Delta_C}
        \ar[ddrrrr]
     &&&&&&&& C \\
    &&&& \ar@{=>}[uu]^i  \\
    &&&& T \ar[uurrrr]_{\id}    }
\]
\end{description}
making the following diagrams commute:
\begin{description}
    \item[{${\bf A1'}$}]
\[
\def\objectstyle{\scriptstyle}
\def\labelstyle{\scriptstyle}
 \xy
    (-45,0)*{(m \ten 1 \ten 1) \circ (m \ten 1) \circ m}="1";
    (-45,-13)*{ (1 \ten m \ten 1)  \circ (m \ten 1) \circ m}="2";
    (45,-13)*{ (1 \ten m \ten 1) \circ ( 1 \ten m ) \circ m}="3";
    (45,0)*{(1 \ten 1 \ten m) \circ (1 \ten m) \circ m }="4";
    (0,0)*{(m \ten m) \circ m}="5";
    (-8,0)*{}="5'";
        {\ar^{a \ten 1 } "1";"2"};
        {\ar^>>>>>>>>>>>>>{(m \ten 1 \ten 1) \circ a} "1";"5"};
        {\ar_{(1 \ten m \ten 1) \circ a} "2";"3"};
        {\ar^{1 \ten a } "3";"4"};
        {\ar^<<<<<<<<<<<<<{(1 \ten 1 \ten m) \circ a} "5";"4"};
 \endxy
\]
    \item[{${\bf A2'}$}]
 \[
 \def\objectstyle{\scriptstyle}
\def\labelstyle{\scriptstyle}
  \xymatrix{
     (1 \ten \id \ten 1) \circ (m \ten 1) \circ m
         \ar[rr]^{(1 \ten \id \ten 1) \circ a}
         \ar[dr]_{r \ten 1 }
     &&  (1 \ten id \ten 1) \circ(1 \ten m) \circ m
         \ar[dl]^{1 \ten \ell}     \\
     & m   }
 \]
    \item[{${\bf A3'}$}]
 \[
\def\objectstyle{\scriptstyle}
\def\labelstyle{\scriptstyle}
 \xy
    (-45,0)*{(\id \ten 1) \circ m}="1";
    (-45,-13)*{1}="2";
    (45,-13)*{(1 \ten \inv) \circ m }="3";
    (45,0)*{(1 \ten \inv \ten 1) \circ (1 \ten m) \circ m }="4";
    (-5,0)*{(1 \ten \inv \ten 1) \circ (m \ten 1) \circ m}="5";
        {\ar^{\ell } "1";"2"};
        {\ar^>>>>>>>>>>>>>{i \ten 1} "1";"5"};
        {\ar_{r^{-1}} "2";"3"};
        {\ar^{1 \ten e } "4";"3"};
        {\ar^{(1 \ten \inv \ten 1) \circ a} "5";"4"};
 \endxy
\]
    \item[{${\bf A4'}$}]
 \[
\def\objectstyle{\scriptstyle}
\def\labelstyle{\scriptstyle}
 \xy
    (-45,0)*{(\inv \ten id) \circ m}="1";
    (-45,-13)*{\inv}="2";
    (45,-13)*{(\id \ten \inv) \circ m }="3";
    (45,0)*{(\inv \ten 1 \ten \inv) \circ (m \ten 1) \circ m}="4";
    (-5,0)*{(\inv \ten 1 \inv) \circ (1 \ten m) \circ m}="5";
        {\ar^{r } "1";"2"};
        {\ar^>>>>>>>>>>>>>{1 \ten i} "1";"5"};
        {\ar_{\ell^{-1}} "2";"3"};
        {\ar^{e \ten 1 } "4";"3"};
        {\ar^{(1 \ten \inv \ten 1) \circ a^{-1}} "5";"4"};
 \endxy
\]
\end{description}
\end{defn}

\begin{prop} \et
A coherent 2-group object in $\Cat$ is a coherent 2-group.
\end{prop}
\textbf{Proof. } We prove this merely by noting that
that the morphisms in $\Cat$ are functors and the
2-morphisms are natural transformations.  With these substitutions
Definition \ref{coherentobject} becomes Definition \ref{defc}.
\qed

We can define a {\bf topological 2-group} to be a coherent
2-group object in $\TopCat$, the 2-category of topological categories.
Similarly, we can define a {\bf Lie 2-group} to be a coherent 2-group
object in $\DiffCat$, the 2-category of smooth categories.  
Here `topological categories' are categories internal to $\Top$,
while `smooth categories' are categories internal to $\Diff$;
concepts of this sort were first studied by Ehresmann 
\cite{Ehresmann,Ehresmann2}
Baez has used strict Lie 2-groups in his work on
categorified gauge theory \cite{Baez}, and it should be interesting
to extend these ideas to general Lie 2-groups.

It seems difficult to internalize the notion of a weak 2-group as we
have just done for coherent 2-groups.  Naively, the definition of a
`weak 2-group object' would require that $\inv$ be a morphism in the ambient
2-category $K$.  In the case $K = \Cat$ this means that $\inv$ is a functor.
However, requiring that $\inv$ be a functor implies the zig-zag
identities, leading us back to the notion of coherent 2-group.

\section{Improvement} \label{improvementsection}

In this section we show that any weak 2-group can be improved to a
coherent one, using the technique of string diagrams \cite{JS0,Street}.
In a strict monoidal category, one can interpret any string diagram as a
morphism in a unique way.  With the help of Mac Lane's coherence theorem
\cite{MacLane} we can also do this in a weak monoidal category.  To do
this, we interpret any string of objects and $1$'s as a tensor product
of objects where all parentheses start in front and all 1's are removed.
Using the associator and left/right unit laws to do any necessary
reparenthesization and introduction or elimination of 1's, any string
diagram then describes a morphism between tensor products of this sort.
The fact that this morphism is unambiguously defined follows from Mac
Lane's coherence theorem.

Let $C$ be a weak 2-group. We will only need string diagrams where
all the strings are labelled by $x$ and $\xb$, where $x$ is some
fixed object of $C$ and $\xb$ is a chosen weak inverse for $x$.
Thus we can omit these labels and just use downwards or upwards
arrows on our strings to distinguish between $x$ and $\xb$. We fix
isomorphisms $i_x \maps$ $\xymatrix@1{1 \ar[r]^<<<<<{\sim} & x
\ten \xb}$ and $e_x \maps$ $\xymatrix@1{\xb \ten x
\ar[r]^<<<<<{\sim} & 1}$. We draw these just as we did in Section
\ref{internalizationsection}; however, they need not satisfy the
zig-zag identities.  Nonetheless, if we write $i_x$ as
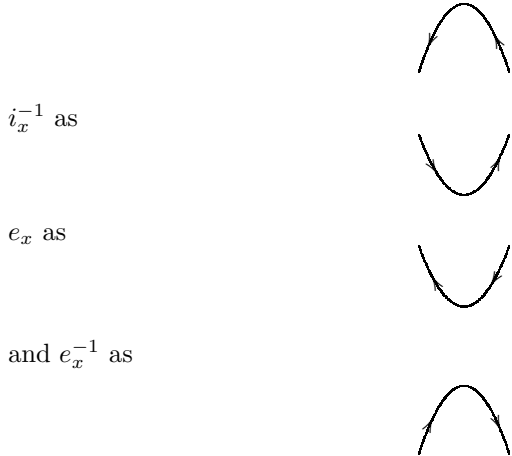
\[
\xy
    (-6,0)*{};
    (6,0)*{};
      **\crv{(0,18)} ?(.16)*\dir{>} ?(.9)*\dir{>};
    (1,11)*{};
\endxy
\]
$i_x^{-1}$ as
\[
 \xy
    (-6,4)*{};(6,4)*{};
      **\crv{(0,-12)}
       ?(.1)*\dir{<}+(2,-1)
       ?(.80)*\dir{<}+(-2,-1);
      (1,-6)*{};
\endxy
\]
$e_x$ as
\[
 \xy
    (-6,4)*{};(6,4)*{};
      **\crv{(0,-12)}
       ?(.20)*\dir{>}+(2,-1)
       ?(.85)*\dir{>}+(-2,-1);
      (1,-6)*{};
\endxy
\]
and $e_x^{-1}$ as
\[
\xy
    (-6,0)*{};
    (6,0)*{};
      **\crv{(0,18)}
      ?(.10)*\dir{<}
      ?(.85)*\dir{<};
    (1,11)*{};
\endxy
\]
we obtain some rules for manipulating string diagrams just
from the fact that these morphisms are inverses of each other.
For instance, the equations
$i_xi_x^{-1} = 1_1$ and $e_x^{-1}e_x = 1_1$
give the rules
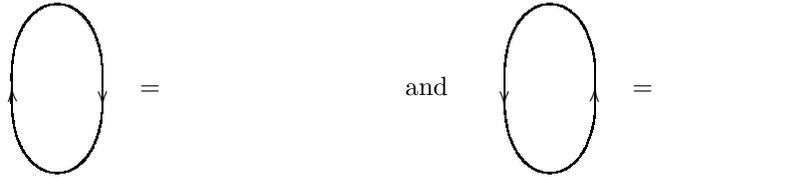
\[
 \xy
    (6,0)*{}="1";
    (-6,0)*{}="2";
      "2";"1" **\crv{(-7,15)&(7,15)}
        ?(0)*\dir{>}  ;
      "1";"2" **\crv{(7,-15) & (-7,-15)}
        ?(.05)*\dir{>} ;
 \endxy
 \quad = \quad
 \qquad \qquad \qquad \qquad {\rm and} \qquad
 \xy
    (6,0)*{}="1";
    (-6,0)*{}="2";
      "1";"2" **\crv{(7,15)&(-7,15)}
        ?(0)*\dir{>}  ;
      "2";"1" **\crv{ (-7,-15)& (7,-15)}
        ?(.05)*\dir{>} ;
 \endxy
 \quad = \quad
 \qquad \qquad .
\]
These rules mean that in a
string diagram, a loop of either form may be removed or inserted
without changing the morphism described by the diagram.
Similarly, the equations $e_x e_x^{-1} = 1_{\bar{x} \tensor x}$
and $i_x^{-1}i_x = 1_{x \tensor \bar{x}}$ give the rules
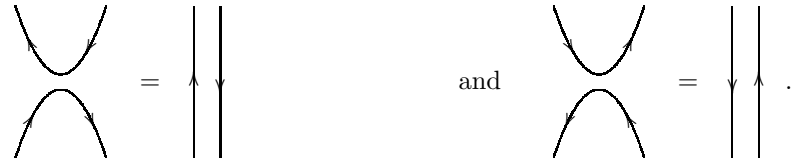
\[
\xy
    (-6,10)*{};
    (6,10)*{};
      **\crv{(0,-8)} ?(.22)*\dir{>} ?(.87)*\dir{>};
    (6,-10)*{};
    (-6,-10)*{};
      **\crv{(0,8)} ?(.20)*\dir{>} ?(.87)*\dir{>};
\endxy
\quad = \quad
 \xy
  (0,-10)*{};(0,10)**\dir{-} ?(.57)*\dir{>};
 \endxy
 \quad
 \xy
  (0,10)*{};(0,-10)**\dir{-} ?(.57)*\dir{>};
 \endxy
 \quad
 \qquad \qquad \qquad \qquad {\rm and} \qquad
 \xy
    (6,10)*{};
    (-6,10)*{};
      **\crv{(0,-8)} ?(.22)*\dir{>} ?(.87)*\dir{>};
    (-6,-10)*{};
    (6,-10)*{};
      **\crv{(0,8)} ?(.20)*\dir{>} ?(.87)*\dir{>};
\endxy
\quad = \quad
 \xy
  (0,10)*{};(0,-10)**\dir{-} ?(.57)*\dir{>};
 \endxy
 \quad
 \xy
  (0,-10)*{};(0,10)**\dir{-} ?(.57)*\dir{>};
 \endxy \quad .
\]
Again, these rules means that in a string diagram we can modify any
portion as above without changing the morphism in question.  We
shall also need another rule, the `horizontal slide':
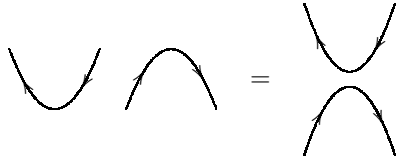
\[
 \xy
    (-6,4)*{};(6,4)*{};
      **\crv{(0,-12)}
       ?(.20)*\dir{>}+(2,-1)
       ?(.85)*\dir{>}+(-2,-1);
\endxy
\quad
 \xy
    (6,-4)*{}; (-6,-4)*{};
      **\crv{(0,12)}
        ?(.20)*\dir{>}+(2,1)
        ?(.85)*\dir{>}+(-2,1);
\endxy
\quad = \quad
 \xy
    (-6,10)*{};
    (6,10)*{};
      **\crv{(0,-8)} ?(.22)*\dir{>} ?(.87)*\dir{>};
    (6,-10)*{};
    (-6,-10)*{};
      **\crv{(0,8)} ?(.20)*\dir{>} ?(.87)*\dir{>};
\endxy
\]
This follows from general results on string diagrams \cite{JS},
but it is easy to prove directly.  First, write down the
corresponding globular diagram:
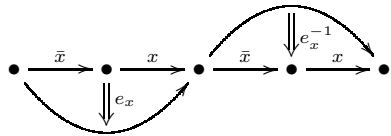
\[
\\
\\
 \xymatrix{
    \bullet
      \ar[r]^{\bar{x}}
      \ar@/_2pc/[rr]_{\quad}_{}="1"
  & \bullet
        \ar@{=>} {};"1"^{e_x}
      \ar[r]^{x}
  & \bullet
      \ar[r]^{\bar{x}}
      \ar@/^2pc/[rr]_{\quad}_{}="1"
  & \bullet
      \ar[r]^{x}
      \ar@{=>}"1" ;{}^<<<<{e_x^{-1}}
  & \bullet
  }
\]
Then, use the bicategory axioms \cite{Benabou} and
Mac Lane's coherence theorem to manipulate this
diagram in the following manner:
\ban
 \xymatrix{
 \def\labelstyle{\scriptstyle}
    \bullet
      \ar[r]^{\xb}
      \ar@/_2pc/[rr]_{\quad}_{}="1"
  & \bullet
        \ar@{=>} {};"1" ^{e_x}
      \ar[r]^{x}
  & \bullet
      \ar[r]^{\xb}
      \ar@/^2pc/[rr]_{\quad}_{}="1"
  & \bullet
      \ar[r]^{x}
      \ar@{=>}"1" ;{}^<<<<{e_x^{-1}}        
  & \bullet
  }
  \quad & = & \quad
 \xymatrix{
    \bullet
      \ar@/^2pc/[rr]^{\xb x}|{\bullet}_{}="0"
      \ar[r]
      \ar@/_2pc/[rr]_{1}_{}="2"           
 &  \bullet
      \ar[r]
      \ar@{=>}"0";{} ^{1_{\xb x}}
      \ar@{=>}{};"2" ^{e_x}
 &  \bullet
      \ar@/^2pc/[rr]^{1}_{}="0"
      \ar[r]
      \ar@/_2pc/[rr]_{\xb x}|{\bullet}_{}="2"
 &  \bullet
      \ar[r]
      \ar@{=>}"0";{} ^<<<{e_x^{-1}}
      \ar@{=>}{};"2" ^{1_{\xb x}}
 &  \bullet
 } \\
  & = & \quad
 \xymatrix{
    \bullet
      \ar@/^2pc/[rr]^{\xb x}|{\bullet}_{}="0"
      \ar@/_2pc/[rr]_{1}_{}="2"
      \ar@{=>}"0";"2" ^{e_{x} \cdot 1_{\xb x}}          
 &&  \bullet
      \ar@/^2pc/[rr]^{1}_{}="0"
      \ar@/_2pc/[rr]_{\xb x}|{\bullet}_{}="2"
      \ar@{=>}"0";"2" ^{1_{\xb x} \cdot e_x^{-1}}
 &&  \bullet
 } \\
 & = &  \quad
 \xymatrix{
     \bullet
      \ar@/^2pc/[rr]^{\xb x}|{\bullet}_{}="0"
      \ar[rr]_{1}="1"
      \ar@{=>}"0";"1" ^{e_x^{-1}}
 &&  \bullet
      \ar[rr]^{1}="1"                                     
      \ar@/_2pc/[rr]_{\xb x}|{\bullet}_{}="2"
      \ar@{=>}"1";"2" ^{e_x}
 &&  \bullet
 } \\
 & = &  \quad
 \xymatrix{
     \bullet
      \ar@/^2pc/[rr]^{\xb x}|{\bullet}_{}="0"
      \ar[rr]_<<<<<{1}^{}="1"
      \ar@/_2pc/[rr]_{1}_{}="2"
      \ar@{=>}"0";"1" ^{e_x^{-1}}
      \ar@{=>}"1";"2" ^{1_1}
 &&  \bullet
      \ar@/^2pc/[rr]^{1}_{}="0"
      \ar[rr]^<<<<<{1}_{}="1"                                     
      \ar@/_2pc/[rr]_{\xb x}|{\bullet}_{}="2"
      \ar@{=>}"0";"1" ^{1_1}
      \ar@{=>}"1";"2" ^{e_x}
 &&  \bullet
 } \\
 & = &  \quad
 \xymatrix{
    \bullet
      \ar@/^2pc/[rr]^{\xb x}_{}="0"
      \ar[r]
      \ar@/_2pc/[rr]_{\xb x}_{}="2"                                       
 &  \bullet
      \ar@{=>}"0";{} ^{e_{x}}
      \ar@{=>}{};"2" ^{e_{x}^{-1}}
      \ar[r]
 &  \bullet
 }
\ean
Rewriting the end result as a string diagram, the result follows.
A similar argument proves another version of the horizontal slide:
\[
 \xy
    (6,4)*{};(-6,4)*{};
      **\crv{(0,-12)}
       ?(.20)*\dir{>}+(2,-1)
       ?(.85)*\dir{>}+(-2,-1);
\endxy
\quad
 \xy
     (-6,-4)*{};(6,-4)*{};
      **\crv{(0,12)}
        ?(.20)*\dir{>}+(2,1)
        ?(.85)*\dir{>}+(-2,1);
\endxy
\quad = \quad
 \xy
    (6,10)*{};
    (-6,10)*{};
      **\crv{(0,-8)} ?(.22)*\dir{>} ?(.87)*\dir{>};
    (-6,-10)*{};
    (6,-10)*{};
      **\crv{(0,8)} ?(.20)*\dir{>} ?(.87)*\dir{>};
\endxy .
\]
With the help of these rules we now prove:

\begin{thm} \label{improve} \et
Given any weak 2-group $C$, it can be improved to a coherent
2-group $\imp(C)$ by equipping each object with an adjoint
equivalence.
\end{thm}

\noindent \textbf{Proof. } First, for each object $x$ we choose a
weak inverse $\xb$ and isomorphisms $i_x \maps 1 \to x \ten \xb$,
$e_x \maps \xb \ten x \to 1$. From this data we construct an
adjoint equivalence $(x,\xb,i'_x,e_x)$ where $i'_x$ is defined as
the following composite morphism:
\[
\xymatrix {1 \ar[r]^{i_x} & x\bar{x} \ar[r]^{ x \cdot
\ell^{-1}_{\bar{x}}} & x(1\bar{x}) \ar[r]^{x \cdot e^{-1}_{x}
\cdot \bar{x}} & x((\bar{x}x)\bar{x}) \ar[r]^{x \cdot
a_{\bar{x},x,\bar{x}}} & x(\bar{x}(x\bar{x}))
\ar[r]^{a^{-1}_{x,\bar{x},x\bar{x}}} & (x\bar{x})(x\bar{x}) \\
\ar[r]^>>>>{i^{-1}_x \cdot (x\bar{x})} & 1(x\bar{x})
\ar[r]^{a^{-1}_{1,x,\bar{x}}} & (1x)\bar{x} \ar[r]^{\ell_x \cdot
\bar{x}} &x\bar{x} . }
\]
where we omit tensor product symbols for brevity.

The above rather cryptic formula for $i'_x$ becomes much clearer if we
use pictures.  If we think of a weak 2-group as a one-object bicategory
and write this formula in globular notation it becomes:
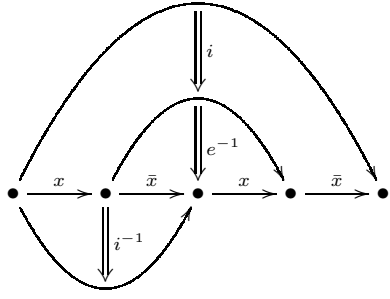
\[
\xymatrix{
     \bullet 
  \ar[r]^x
    \ar@/^6pc/[rrrr]^{\quad}_{}="0"
    \ar@/_3pc/[rr]^{\quad}_{}="3"
 &  \bullet 
  \ar[r]^{\bar{x}}
      \ar@/^3pc/[rr]^{\quad}_{}="1"
           \ar@{=>} {} ;"3"  ^{i^{-1}}
 &  \bullet 
  \ar[r]^x
           \ar@{=>}"1"; {} ^{e^{-1}}
 &  \bullet 
  \ar[r]^{\bar{x}}
           \ar@{=>}"0"; "1" ^{i}
 &  \bullet 
    }
\]
If we write it as a string diagram it looks even simpler:
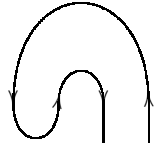
\[
 \xy
    (-6,-0)*{}="1";
    (0,0)*{}="2";
    (6,0)*{}="3";
    (6,-5)*{}="3'";
    (12,0)*{}="4";
    (12,-5)*{}="4'";
      "4";"1" **\crv{(12,18)& (-6,18)}
        ?(.03)*\dir{>}  ?(1)*\dir{>};
      "1";"2" **\crv{(-6,-6)&(0,-6)};
      "2";"3" **\crv{(0,6)&(6,6)};
        ?(.08)*\dir{>} ;
      "3";"3'" **\dir{-}
        ?(0)*\dir{>};
      "4";"4'" **\dir{-};
 \endxy
\]

Now let us show that
$(x, \bar{x}, i'_x, e_x)$ satisfies the zig-zag identities.
Recall that these identities say that the following diagrams commute:
\[
 \xymatrix{
     1 \tensor x \ar[r]^{i'_x \tensor 1} \ar[d]_{\ell_{x}}
     & (x \tensor \bar{x}) \tensor x \ar[r]^{a_{x, \bar{x},x} }
     & x \tensor ( \bar{x} \tensor x) \ar[d]^{1 \tensor e_x} \\
     x \ar[rr]_{r^{-1}_x}
      && x \tensor 1  }
\qquad
 \xymatrix{
    \bar{x} \ten 1 \ar[r]^{1 \ten i'_x} \ar[d]_{r_{\bar{x}}}
     & \bar{x} \ten (x \ten \bar{x}) \ar[r]^{a^{-1}_{\bar{x},x, \bar{x}}}
     & (\bar{x} \ten x) \ten \bar{x} \ar[d]^{e_x \ten 1} \\
    \bar{x} \ar[rr]_{\ell^{-1}_{\bar{x}}}
     && 1 \ten x   } \quad .
\]

Utilizing the observation made regarding Mac Lane's coherence
theorem we can express the zig-zag identities in globular notation
as follows:
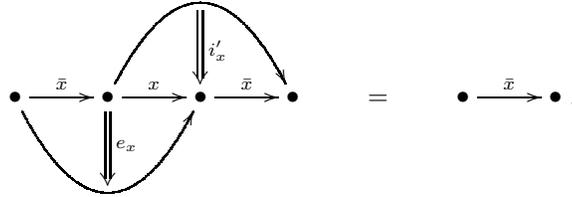
\[
\xymatrix{
   \bullet
     \ar[r]^x
     \ar@/^3pc/[rr]^{\quad}_{}="0"
 & \bullet
     \ar[r]^{\bar{x}}
     \ar@/_3pc/[rr]^{\quad}_{}="2"
     \ar@{=>}"0"; {} ^{i'_x}
 & \bullet
     \ar[r]^x
     \ar@{=>}{}; "2" ^{e_x}
 & \bullet }
 \qquad = \qquad
\xymatrix{ \bullet \ar[r]^x & \bullet } ,
 \]
\[
\xymatrix{
   \bullet
     \ar[r]^{\bar{x}}
     \ar@/_3pc/[rr]^{\quad}_{}="0"
 & \bullet
     \ar[r]^x
     \ar@/^3pc/[rr]^{\quad}_{}="2"
     \ar@{=>};"0" ^{e_x}
 & \bullet
     \ar[r]^{\bar{x}}
     \ar@{=>}"2";{} ^{i'_x}
 & \bullet }
\qquad = \qquad \xymatrix{ \bullet \ar[r]^{\bar{x}} & \bullet } .
\]
If we express $i'_x$ in terms of $i_x$ and $e_x$, these
equations become
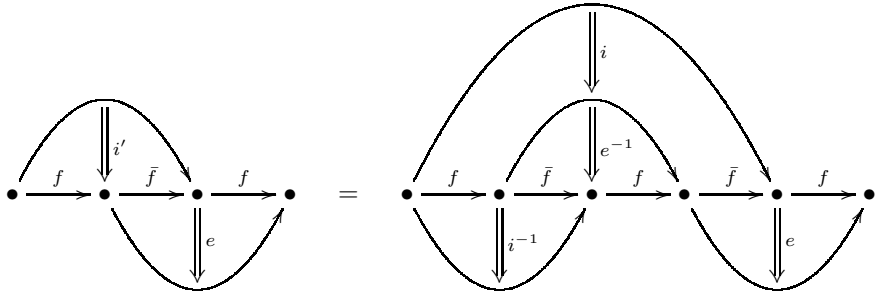
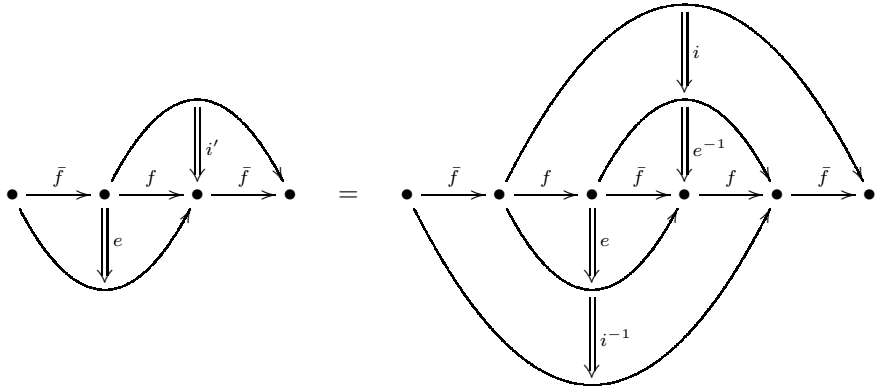
\[
\xymatrix{
   \bullet
     \ar[r]^f
     \ar@/^3pc/[rr]^{\quad}_{}="0"
 & \bullet
     \ar[r]^{\bar{f}}
     \ar@/_3pc/[rr]^{\quad}_{}="2"
     \ar@{=>}"0"; {} ^{i'}
 & \bullet
     \ar[r]^f
     \ar@{=>}{}; "2" ^{e}
 & \bullet }
\quad = \quad
 \xymatrix{
     \bullet 
  \ar[r]^f
    \ar@/^6pc/[rrrr]^{\quad}_{}="0"
    \ar@/_3pc/[rr]^{\quad}_{}="3"
 &  \bullet 
  \ar[r]^{\bar{f}}
      \ar@/^3pc/[rr]^{\quad}_{}="1"
           \ar@{=>} {} ;"3"  ^{i^{-1}}
 &  \bullet 
  \ar[r]^f
           \ar@{=>}"1"; {} ^{e^{-1}}
 &  \bullet 
  \ar[r]^{\bar{f}}
      \ar@/_3pc/[rr]^{\quad}_{}="3"
           \ar@{=>}"0"; "1" ^{i}
 &  \bullet 
  \ar[r]^f
          \ar@{=>} {} ;"3"  ^{e}
 & \bullet   }
\\
\\
\]
and
\[
\xymatrix{
   \bullet
     \ar[r]^{\bar{f}}
     \ar@/_3pc/[rr]^{\quad}_{}="0"
 & \bullet
     \ar[r]^f
     \ar@/^3pc/[rr]^{\quad}_{}="2"
     \ar@{=>};"0" ^{e}
 & \bullet
     \ar[r]^{\bar{f}}
     \ar@{=>}"2";{} ^{i'}
 & \bullet }
\quad = \quad
 \xymatrix{
     \bullet 
  \ar[r]^{\bar{f}}
        \ar@/_6pc/[rrrr]^{\quad}_{}="2"
 &  \bullet 
  \ar@/^6pc/[rrrr]^{\quad}_{}="0"
  \ar[r]^f
  \ar@/_3pc/[rr]^{\quad}_{}="3"
 &  \bullet 
  \ar[r]^{\bar{f}}
  \ar@/^3pc/[rr]^{\quad}_{}="1"
  \ar@{=>}{}; "3" ^{e}
 &  {\bullet} 
  \ar[r]^f
   \ar@{=>}"0";"1"^{i}
   \ar@{=>}"1"; {} ^{e^{-1}}
   \ar@{=>}"3"; "2" ^{i^{-1}}
 &  \bullet 
  \ar[r]^{\bar{f}}
 & \bullet   }
\]

To verify these two equations we use string diagrams.  In the
calculations that follow, we denote an application of the
`horizontal slide' rule by a dashed line
connecting the appropriate zig and zag.  Dotted lines connecting two
parallel strings will indicate an application of the rules $e_x e_x^{-1}
= 1_{\bar{x} \tensor x}$ or $i_x^{-1}i_x = 1_{x \tensor \bar{x}}$.
Furthermore, the rules $i_xi_x^{-1} = 1_1$ and $e_x^{-1}e_x = 1_1$ allow
us to remove a closed loop any time one appears.  The first equation
can be proved as follows:
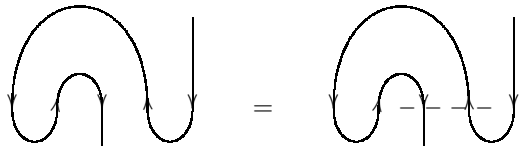
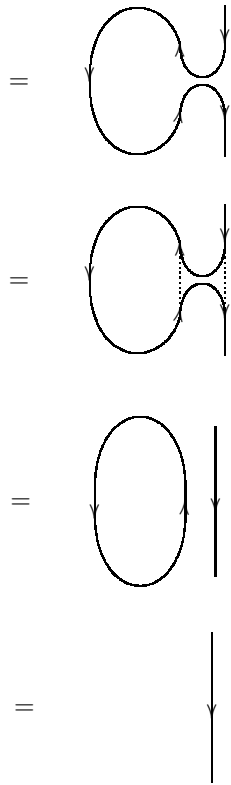
\[
 \xy
    (-6,-0)*{}="1";
    (0,0)  *{}="2";
    (6,0)  *{}="3";
    (6,-5) *{}="3'";
    (12,0) *{}="4";
    (18,0) *={}="5";
    (18,12)*={}="5'";
      "5";"4" **\crv{(18,-6)& (12,-6)};
      "5'";"5" **\dir{-} ?(1)*\dir{>};
      "4";"1" **\crv{(12,18)& (-6,18)};
        ?(.03)*\dir{>}  ?(1)*\dir{>};
      "1";"2" **\crv{(-6,-6)&(0,-6)};
      "2";"3" **\crv{(0,6)&(6,6)};
        ?(.08)*\dir{>} ;
      "3";"3'" **\dir{-}
        ?(0)*\dir{>};
 \endxy
\qquad = \qquad
 \xy
    (-6,-0)*{}="1";
    (0,0)*{}="2";
    (6,0)*{}="3";
    (6,-5)*{}="3'";
    (12,0)*{}="4";
    (18,0)*={}="5";
    (18,12)*={}="5'";
      "5";"4" **\crv{(18,-6)& (12,-6)};
      "5'";"5" **\dir{-} ?(1)*\dir{>};
      "4";"1" **\crv{(12,18)& (-6,18)}
        ?(.03)*\dir{>}  ?(1)*\dir{>};
      "1";"2" **\crv{(-6,-6)&(0,-6)};
      "2";"3" **\crv{(0,6)&(6,6)};
        ?(.08)*\dir{>} ;
      "3";"3'" **\dir{-}
        ?(0)*\dir{>} ;
    (3,0)*{}="A";
    (15,0)*{}="B";
        "A";"B" **\dir{--};
  \endxy
\]
\[
 \xy
    (-6,-0)*{}="1";
    (0,0)  *{}="2";
    (6,0)  *{}="3";
    (6,-5) *{}="3'";
    (12,0) *{}="4";
    (18,0) *={}="5";
    (18,12)*={}="5'";
 \endxy
\qquad = \qquad \xy
    (-6,-0)*{}="1";
    (6,-5)  *{}="2";
    (12,-5)  *{}="3";
    (12,-10) *{}="3'";
    (6,5) *{}="4";
    (12,5) *={}="5";
    (12,10)*={}="5'";
    (18,0) *={}="";
      "5";"4" **\crv{(12,-1)& (6,-1)}
            ?(1)*\dir{>};
      "5'";"5" **\dir{-} ?(1)*\dir{>};
      "4";"1" **\crv{(4,12)& (-6,12)};
          ?(1)*\dir{>};
      "1";"2" **\crv{(-6,-12)&(4,-12)};
      "2";"3" **\crv{(6,1)&(12,1)};
        ?(.08)*\dir{>} ;
      "3";"3'" **\dir{-}
        ?(0)*\dir{>};
 \endxy
\]
\[
 \xy
    (-6,-0)*{}="";
    (18,0) *={}="";
 \endxy
\qquad = \qquad \xy
    (-6,-0)*{}="1";
    (6,-5)  *{}="2";
    (12,-5)  *{}="3";
    (12,-10) *{}="3'";
    (6,5) *{}="4";
    (12,5) *={}="5";
    (12,10)*={}="5'";
    (18,0) *={}="";
      "5";"4" **\crv{(12,-1)& (6,-1)}
            ?(1)*\dir{>};
      "5'";"5" **\dir{-} ?(1)*\dir{>};
      "4";"1" **\crv{(4,12)& (-6,12)};
          ?(1)*\dir{>};
      "1";"2" **\crv{(-6,-12)&(4,-12)};
      "2";"3" **\crv{(6,1)&(12,1)};
        ?(.08)*\dir{>} ;
      "3";"3'" **\dir{-}
        ?(0)*\dir{>};
      "2";"4" **\dir{.};
      "3";"5" **\dir{.};
 \endxy
\]
\[
 \xy
    (-6,-0)*{}="";
    (18,0) *={}="";
 \endxy
\qquad = \qquad
 \xy
    (6,0)*{}="1";
    (-6,0)*{}="2";
      "1";"2" **\crv{(7,15)&(-7,15)}
        ?(0)*\dir{>}  ;
      "2";"1" **\crv{ (-7,-15)& (7,-15)}
        ?(.05)*\dir{>} ;
 \endxy
\quad
 \xy
  (0,10)*{};(0,-10)**\dir{-} ?(.57)*\dir{>};
 \endxy \qquad
\]
\[
 \xy
    (-6,-0)*{}="1";
    (18,0) *={}="5";
 \endxy
\qquad = \qquad
 \xy
    (6,0)*{}="";
    (-6,0)*{}="";
 \endxy
\quad
 \xy
  (0,10)*{};(0,-10)**\dir{-} ?(.57)*\dir{>};
 \endxy \qquad
\]

The proof of the second equation is accomplished in
a similar manner:
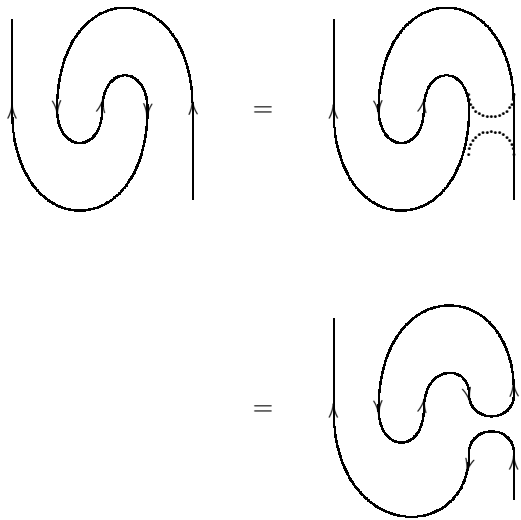
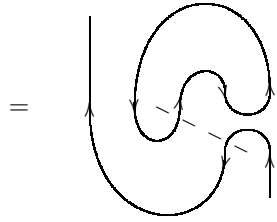
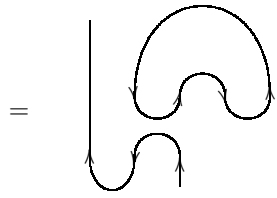
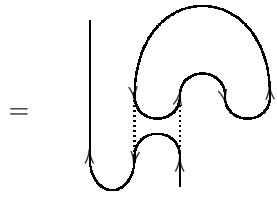
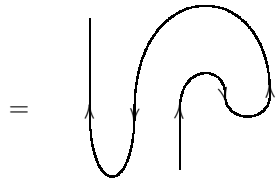
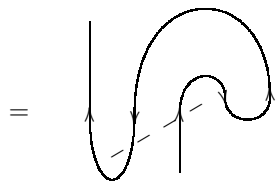
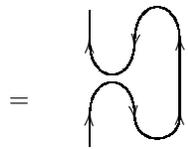
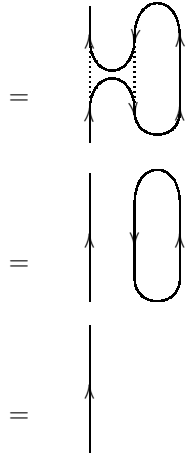
\[
 \xy
    (-6,-0)*{}="1";
    (0,0)*{}="2";
    (6,0)*{}="3";
    (12,0)*{}="4";
    (12,-12)*{}="4'";
    (-12,12)*{}="0'";
    (-12,0)*{}="0";
      "4";"1" **\crv{(12,18)& (-6,18)};
        ?(.02)*\dir{>}  ;
      "1";"2" **\crv{(-6,-6) & (0,-6)}
        ?(.04)*\dir{>};
      "2";"3" **\crv{(0,6) & (6,6)};
        ?(.08)*\dir{>} ;
      "4'";"4" **\dir{-};
      "3";"0" **\crv{(6,-18)&(-12,-18)}
        ?(.02)*\dir{>} ;
      "0";"0'" **\dir{-} ?(.05)*\dir{>};
 \endxy
\qquad = \qquad
 \xy
    (-6,-0)*{}="1";
    (0,0)*{}="2";
    (6,0)*{}="3";
    (12,0)*{}="4";
    (12,-12)*{}="4'";
    (-12,12)*{}="0'";
    (-12,0)*{}="0";
    (6,2)*{}="A";
    (6,-6)*{}="A'";
    (12,2)="B";
    (12,-6)="B'";
      "4";"1" **\crv{(12,18)& (-6,18)};
      "1";"2" **\crv{(-6,-6) & (0,-6)}
        ?(.04)*\dir{>};
      "2";"3" **\crv{(0,6) & (6,6)};
        ?(.08)*\dir{>} ;
      "4'";"4" **\dir{-};
      "3";"0" **\crv{(6,-18)&(-12,-18)};
      "0";"0'" **\dir{-} ?(.05)*\dir{>};
      "A";"B" **\crv{~*=<2pt>{.}(6,-2)&(12,-2)};
      "A'";"B'" **\crv{~*=<2pt>{.}(6,-2)&(12,-2)};
 \endxy
\]
\[
 \xy
    (12,0)*{};
    (-12,0)*{};
 \endxy
\qquad = \qquad
 \xy
    (-6,-0)*{}="1";
    (0,0)*{}="2";
    (6,0)*{}="3";
    (12,0)*{}="4";
    (12,-12)*{}="4'";
    (-12,12)*{}="0'";
    (-12,0)*{}="0";
    (6,2)*{}="A";
    (6,-6)*{}="A'";
    (12,2)="B";
    (12,-6)="B'";
      "B";"1" **\crv{(12,18)& (-6,18)}
        ?(.03)*\dir{>};
      "1";"2" **\crv{(-6,-6) & (0,-6)}
        ?(.04)*\dir{>};
      "2";"A" **\crv{(0,6) & (6,6)};
        ?(.08)*\dir{>} ;
      "4'";"B'" **\dir{-}
         ?(.96)*\dir{>};
      "A'";"0" **\crv{(6,-18)&(-12,-18)}
        ?(.07)*\dir{>};
      "0";"0'" **\dir{-} ?(.05)*\dir{>};
      "A";"B" **\crv{(6,-2)&(12,-2)}
        ?(.08)*\dir{>} ;
      "B'";"A'" **\crv{(12,-2)&(6,-2)};
 \endxy
\]
\[
 \xy
    (12,0)*{};
    (-12,0)*{};
 \endxy
\qquad = \qquad
 \xy
    (-6,-0)*{}="1";
    (0,0)*{}="2";
    (6,0)*{}="3";
    (12,0)*{}="4";
    (12,-12)*{}="4'";
    (-12,12)*{}="0'";
    (-12,0)*{}="0";
    (6,2)*{}="A";
    (6,-6)*{}="A'";
    (12,2)="B";
    (12,-6)="B'";
      "B";"1" **\crv{(12,18)& (-6,18)};
        ?(.03)*\dir{>};
      "1";"2" **\crv{(-6,-6) & (0,-6)};
        ?(.04)*\dir{>} ;
      "2";"A" **\crv{(0,6) & (6,6)};
        ?(.09)*\dir{>};
      "4'";"B'" **\dir{-};
         ?(.96)*\dir{>};
      "A'";"0" **\crv{(6,-18)&(-12,-18)};
        ?(.07)*\dir{>} ;
      "0";"0'" **\dir{-} ?(.05)*\dir{>};
      "A";"B" **\crv{(6,-2)&(12,-2)};
        ?(.08)*\dir{>};
      "B'";"A'" **\crv{(12,-2)&(6,-2)};
    (-3,0)*{}="X";
    (9,-6)*{}="XX";
    "X";"XX" **\dir{--};
 \endxy
\]
\[
 \xy
    (12,0)*{};
    (-12,0)*{};
 \endxy
\qquad = \qquad
 \xy
    (-6,2)*{}="1";
    (0,2)*{}="2";
    (6,0)*{}="3";
    (12,0)*{}="4";
    (0,-10)*{}="4'";
    (-12,12)*{}="0'";
    (-12,-6)*{}="0";
    (6,2)*{}="A";
    (-6,-6)*{}="A'";
    (12,2)="B";
    (0,-6)="B'";
      "B";"1" **\crv{(12,18)& (-6,18)}
        ?(.03)*\dir{>};
      "1";"2" **\crv{(-6,-2) & (0,-2)}
        ?(.04)*\dir{>};
      "2";"A" **\crv{(0,6) & (6,6)};
        ?(.05)*\dir{>} ;
      "4'";"B'" **\dir{-}
         ?(.96)*\dir{>};
      "A'";"0" **\crv{(-6,-12)&(-12,-12)}
        ?(.07)*\dir{>};
      "0";"0'" **\dir{-} ?(.05)*\dir{>};
      "A";"B" **\crv{(6,-2)&(12,-2)}
        ?(.08)*\dir{>} ;
      "B'";"A'" **\crv{(0,-2)&(-6,-2)};
 \endxy
\]
\[
 \xy
    (12,0)*{};
    (-12,0)*{};
 \endxy
\qquad = \qquad
 \xy
    (-6,2)*{}="1";
    (0,2)*{}="2";
    (6,0)*{}="3";
    (12,0)*{}="4";
    (0,-10)*{}="4'";
    (-12,12)*{}="0'";
    (-12,-6)*{}="0";
    (6,2)*{}="A";
    (-6,-6)*{}="A'";
    (12,2)="B";
    (0,-6)="B'";
      "B";"1" **\crv{(12,18)& (-6,18)}
        ?(.03)*\dir{>};
      "1";"2" **\crv{(-6,-2) & (0,-2)}
        ?(.04)*\dir{>};
      "2";"A" **\crv{(0,6) & (6,6)};
        ?(.05)*\dir{>} ;
      "4'";"B'" **\dir{-}
         ?(.96)*\dir{>};
      "A'";"0" **\crv{(-6,-12)&(-12,-12)}
        ?(.07)*\dir{>};
      "0";"0'" **\dir{-} ?(.05)*\dir{>};
      "A";"B" **\crv{(6,-2)&(12,-2)}
        ?(.08)*\dir{>} ;
      "B'";"A'" **\crv{(0,-2)&(-6,-2)};
      "1";"A'" **\dir{.};
      "2";"B'" **\dir{.};
 \endxy
\]
\[
 \xy
    (-12,-0)*{};
    (12,0)*{};
 \endxy
\qquad = \qquad
 \xy
    (-6,-0)*{}="1";
    (0,0)*{}="2";
    (0,-8)*{}="2'";
    (6,0)*{}="3";
    (-12,12)*{}="0'";
    (-12,0)*{}="0";
    (6,2)*{}="A";
    (12,2)="B";
        "0";"0'" **\dir{-} ?(.05)*\dir{>};
        "1";"0" **\crv{(-6,-12) &(-12,-12)}
            ?(.05)*\dir{>};
        "B";"1" **\crv{(12,18)& (-6,18)};
          ?(.03)*\dir{>};
        "A";"B" **\crv{(6,-2)&(12,-2)};
          ?(.08)*\dir{>};
        "2";"A" **\crv{(0,6) & (6,6)};
          ?(.03)*\dir{>};
        "2'";"2" **\dir{-};
 \endxy
\]
\[
 \xy
    (-12,-0)*{};
    (12,0)*{};
 \endxy
\qquad = \qquad
 \xy
    (-6,-0)*{}="1";
    (0,0)*{}="2";
    (0,-8)*{}="2'";
    (6,0)*{}="3";
    (-12,12)*{}="0'";
    (-12,0)*{}="0";
    (6,2)*{}="A";
    (12,2)="B";
        "0";"0'" **\dir{-} ?(.05)*\dir{>};
        "1";"0" **\crv{(-6,-12) &(-12,-12)}
            ?(.05)*\dir{>};
        "B";"1" **\crv{(12,18)& (-6,18)};
          ?(.03)*\dir{>};
        "A";"B" **\crv{(6,-2)&(12,-2)};
          ?(.08)*\dir{>};
        "2";"A" **\crv{(0,6) & (6,6)};
          ?(.03)*\dir{>};
        "2'";"2" **\dir{-};
    (-9,-6)*{}="X";
    (3,1)*{}="XX";
    "X";"XX" **\dir{--};
 \endxy
\]
\[
 \xy
    (-12,-0)*{};
    (12,0)*{};
 \endxy
\qquad = \qquad
 \xy
    (-12,12)*{}="0'";
    (-12,8)*{}="0";
    (-6,8)*{}="1";
    (-12,-2)*{}="2";
    (-12,-6)*{}="2'";
    (6,-2)*{}="3";
    (-6,-2)*{}="A";
    (0,-2)*{}="B";
    (0,8)*{}="B'";
    (12,0)*{}="";
        "0";"0'" **\dir{-} ?(.05)*\dir{>};
        "1";"0" **\crv{(-6,2) &(-12,2)}
            ?(.05)*\dir{>};
        "B'";"1" **\crv{(0,14)& (-6,14)};
          ?(.03)*\dir{>};
        "A";"B" **\crv{(-6,-6)&(0,-6)};
          ?(.04)*\dir{>};
        "2";"A" **\crv{(-12,4) & (-6,4)};
          ?(.03)*\dir{>};
        "2'";"2" **\dir{-};
        "B";"B'" **\dir{-} ?(.06)*\dir{>};;
 \endxy
\]
\[
 \xy
    (-12,-0)*{};
    (12,0)*{};
 \endxy
\qquad = \qquad
 \xy
    (-12,12)*{}="0'";
    (-12,8)*{}="0";
    (-6,8)*{}="1";
    (-12,-2)*{}="2";
    (-12,-6)*{}="2'";
    (6,-2)*{}="3";
    (-6,-2)*{}="A";
    (0,-2)*{}="B";
    (0,8)*{}="B'";
    (12,0)*{}="";
        "0";"0'" **\dir{-} ?(.05)*\dir{>};
        "1";"0" **\crv{(-6,2) &(-12,2)}
            ?(.05)*\dir{>};
        "B'";"1" **\crv{(0,14)& (-6,14)};
          ?(.03)*\dir{>};
        "A";"B" **\crv{(-6,-6)&(0,-6)};
          ?(.04)*\dir{>};
        "2";"A" **\crv{(-12,4) & (-6,4)};
          ?(.03)*\dir{>};
        "2'";"2" **\dir{-};
        "B";"B'" **\dir{-} ?(.06)*\dir{>};
        "0";"2" **\dir{.};
        "1";"A" **\dir{.};
 \endxy
\]
\[
 \xy
    (-12,-0)*{};
    (12,0)*{};
 \endxy
\qquad = \qquad
 \xy
    (-12,12)*{}="0'";
    (-6,8)*{}="1";
    (-12,-5)*{}="2'";
    (6,-2)*{}="3";
    (-6,-2)*{}="A";
    (0,-2)*{}="B";
    (0,8)*{}="B'";
    (12,0)*{}="";
        "B'";"1" **\crv{(0,14)& (-6,14)};
        "A";"B" **\crv{(-6,-6)&(0,-6)};
        "B";"B'" **\dir{-} ?(.60)*\dir{>};
        "2'";"0'" **\dir{-} ?(.53)*\dir{>};
        "1";"A" **\dir{-};  ?(.57)*\dir{>};
 \endxy
\]
\[
 \xy
    (-12,-0)*{};
    (12,0)*{};
 \endxy
\qquad = \qquad
 \xy
    (-12,12)*{}="0'";
    (-12,-5)*{}="2'";
    (12,0)*{}="";
        "2'";"0'" **\dir{-} ?(.53)*\dir{>};
 \endxy
\]
\hskip 30em \qed

We can now make this `improvement' process into a 2-functor
$\imp \maps \wg \to \cg$:

\begin{thm} \et
There exist a 2-functor $ \imp \maps \wg \to \cg$ which sends any
object $C \in \wg$ to $\imp(C) \in \cg$ and acts as the identity
on morphisms and 2-morphisms.
\end{thm}

\textbf{Proof. } The proof of this theorem is a trivial consequence of
Theorem \ref{improve}.  Obviously all domains, codomains, identities
and composites are preserved, since the 1-morphisms and 2-morphisms are
unchanged as a result of Definitions \ref{coherenthomo} and
\ref{coherent2homo}. \qed

On the other hand, there is also a forgetful 2-functor ${\rm F} \maps
\cg \to \wg$, which forgets the extra structure on objects
and acts as the identity on morphisms and 2-morphisms.  It is
easy to see that improvement followed by this forgetful 2-functor
acts as the identity of $\wg$.  On the other hand, applying $F$ and
then $\imp$ to a coherent 2-group $C$ amounts to forgetting the
choice of adjoint equivalence for each object $x \in C$ and then
making a new such choice.  We obtain a new coherent 2-group $C'$,
but it has the same underlying weak monoidal category, so the
identity functor $1_C \maps C \to C'$ is a coherent 2-group
homomorphism from $C$ to $C'$.   It should thus not be surprising
that:

\begin{thm} \et
The 2-functors $ \imp \maps \wg \to \cg$, ${\rm F} \maps \cg \to
\wg$ extend to define a biequivalence between the 2-categories
$\wg$ and $\cg$.
\end{thm}

We omit the proof, but refer the reader to B\'enabou \cite{Benabou} and
Street \cite{Street} for a discussion of the concept of biequivalence.
The upshot is that we can use either weak or coherent 2-groups,
whichever happens to be more convenient at the time, and freely
translate results between the two formalisms.

To conclude, let us summarize {\it why} weak and coherent 2-groups are
not really so different.  At first, the choice of a specified adjoint
equivalence for each object seems like a substantial extra structure to
put on a weak 2-group.  However, Theorem \ref{improve} shows that we can
always succeed in putting this extra structure on any weak 2-group.
Furthermore, while there are many ways to equip a weak 2-group with this
extra structure, there is `essentially' just one way, since the remarks
at the end of Section \ref{coherent} show that this structure is {\it
automatically preserved up to coherent isomorphism} by any homomorphism
of weak 2-groups.

\subsection*{Acknowledgements}  I thank John Baez
and James Dolan for very helpful discussions and correspondence,
and also Miguel Carri\'on \'Alvarez for
his assistance in making string diagrams.


\begin{thebibliography}{7}

\bibitem{Baez} J.\ Baez, Higher Yang--Mills theory, available at
hep-th/0206130.

\bibitem{Benabou} J. B\'{e}nabou, {\sl Introduction to Bicategories},
Lecture Notes in Mathematics {\bf 47}, Springer, New York, 1967, pp. 1-77.

\bibitem{Brown} R.\ Brown, Groupoids and crossed objects in algebraic
topology, {\sl Homology, Homotopy and Applications} {\bf 1} (1999), 1--78.
Available at
\texttt{http://www.math.rutgers.edu/hha/volumes/1999/volume1-1.htm}

\bibitem{BS} R.\ Brown and C.\ B.\ Spencer, $\cal G$-groupoids,
crossed modules, and the classifying space of a topological group,
{\sl Proc.\ Kon.\ Akad.\ v.\ Wet.} {\bf 79} (1976), 296--302.

\bibitem{Ehresmann} C.\ Ehresmann, Cat\'egories topologiques et cat\'egories
diff\'erentiables, in {\sl Colloque G\'eom.\ Diff.\ Globale}, 1959,
Centre Belge Rech.\ Math., Louvain, pp.\ 137--150.  Reprinted in
Charles Ehresmann, {\sl Oeuvres Completes Et Commentees},
Cahiers de Topologie et Geometrie Differentielle, Amien, 1980--1984.

\bibitem{Ehresmann2}
C.\ Ehresmann, Introduction to the theory of structured
categories, Tech.\ Report 10, Univ.\ of Kansas at Lawrence, 1966.
Reprinted in Charles Ehresmann, {\sl Oeuvres Completes Et Commentees},
Cahiers de Topologie et Geometrie Differentielle, Amien, 1980--1984.

\bibitem{moncat} S.\ Eilenberg and G.\ M.\ Kelly, Closed categories,
\textit{Proceedings of the Conference on Categorical Algebra (La Jolla,
Calif., 1965)}, Springer, New York, 1966, pp.\ 421--562.

\bibitem{FB} M.\ Forrester-Barker, Group objects and internal
categories, available as math.CT/0212065.

\bibitem{HKK} K.\ A.\ Hardie, K.\ H.\ Kamps and R.\ W.\ Kieboom,
A homotopy bigroupoid of a topological space, {\sl Appl.\ Cat.\
Str.\ } {\bf 9} (2001) 311--327.

\bibitem{JS0} A.\ Joyal and R.\ Street, The geometry of tensor
calculus, I, {\sl Adv.\ Math.\ } {\bf 88} (1991), 55---112.

\bibitem{JS} A.\ Joyal and R.\ Street, Braided tensor categories,
{\sl Adv.\ Math.\ } {\bf 102} (1993), 20--78.

\bibitem{Kelly} G.\ M.\ Kelly, On Mac Lane's conditions for coherence
of natural associativities, commutativities, etc., {\sl J.\ Algebra}
{\bf 4} (1967), 397--402.

\bibitem{Laplaza} M.\ L.\ Laplaza, Coherence for categories
with group structure: an alternative approach, {\sl J.\ Algebra }
{\bf 84} (1983),  305--323.

\bibitem{MacLane} S. MacLane, Natural associativity and
commutativity, {\sl Rice Univ.\ Stud.} {\bf 49} (1963), 28--46.

\bibitem{Street} R.\ Street, Low-dimensional topology and
higher-order categories, {\sl Proceedings of CT95}, Halifax, July 9-15 1995

\bibitem{Street2} R.\ Street, Categorical structures, in
{\sl Handbook of Algebra}, Vol.\ 1, ed.\ M.\ Hazewinkel, Elsevier
North-Holland, 1995, pp.\ 529--577.

\bibitem{Rivano} N.\ Saavedra Rivano, {\sl Cat\'{e}gories Tannakiennes},
Lecture Notes in Mathematics {\bf 265}, Springer, New York, 1972.

\bibitem{Ulbrich} K.-H.\ Ulbrich, Koh\"arenz in Kategorien mit
Gruppenstruktur, {\sl J.\ Algebra} {\bf 72} (1981), 279--295.

\bibitem{Yetter} D.\ Yetter, TQFTs from homotopy 2-types,
{\sl J.\ Knot Theory Ramifications} {\bf 2} (1993), 113--123.

\end{thebibliography}
\end{document}